\documentclass{ws-ijitdm}
\usepackage{tabularx} 
\usepackage{multicol}
\usepackage{multirow}
\usepackage{float}
\usepackage{hyperref}
\hypersetup{colorlinks, allcolors=black}
\usepackage{graphicx} 
\makeatletter
\def\footnoterule{\kern-3\p@
  \hrule \@width 2in \kern 2.6\p@} 
  
\makeatother
\begin{document}
\renewcommand{\footnotesize}{\fontsize{9pt}{11pt}\selectfont}
\title{Virtual Gap Analysis procedures for Multi-Criteria Decision-Making and Efficiency Analysis Problems}

\author{Fuh-Hwa Franklin Liu \footnotemark{}\footnotetext{Corresponding author; fliu@nycu.edu.tw} }
\address{Professor Emeritus, Department of Industrial Engineering and Management,\\  National Yang Ming Chiao Tung University, Taiwan 300, Republic of China}
\author{Su-Chuan Shih \footnotemark{}\footnotetext{scshih@gm.pu.edu.tw}}
\address{Associate Professor, Department of Business Administration,\\ Providence University, Taiwan 433, Republic of China}

\maketitle
\begin{history}
\received{7 June 2024}
\revised{22 January 2025}
\accepted{24 June 2025}
\end{history}
\begin{abstract} 
Existing multi-criteria decision-making (MCDM) methods often face challenges when evaluating a large number of alternatives, leading to skewed results in selecting the optimal choice. Similarly, conventional efficiency analysis (EA) methods, such as Data Envelopment Analysis (DEA) and Stochastic Frontier Analysis (SFA), often yield incomplete solutions due to their reliance on theoretical assumptions. To address these limitations, we propose a novel EA method that integrates Virtual Gap Analysis (VGA) models to evaluate the performance of each decision-making unit (DMU) in relation to others based on best practices. Unlike DEA and SFA, our VGA models are linear programming-based, assumption-free, and capable of delivering robust and reliable solutions. The proposed method enables each DMU to identify achievable benchmarks for inputs and outputs. Based on the estimated virtual gaps, DMUs are classified as inefficient (with scores below one) or efficient (with scores of one or higher). Additionally, our new MCDM method incorporates existing MCDM techniques to analyze the few identified efficient DMUs, significantly reducing the effort required to select the best DMU.

\keywords{Decision Support System; Multiple Criteria Decision Making; Efficiency Analysis.} 
\noindent \textit{MSC}: 90B50; 90C29; 90C08; 91A80; 91B06.
\end{abstract}
 \renewcommand{\thefootnote}{\roman{footnote}}
\section{Research Background and Objectives} \label{sec:1}
\subsection{The Traditional EA and MCDM Methods} \label{sec:1.1}

Table \ref{table:1} presents a numerical example involving six firms—K, A, B, D, G, and H—that consume resources \(X_1\) and \(X_2\) to produce products \(Y_1\) and \(Y_2\) in varying amounts. The resources and products have different measurement units, and for some firms, market prices per unit of the resources and products are unavailable. Interactions between the resources and products are disregarded. The observed dataset consists exclusively of positive, continuous real numbers. 

The observed dataset can be described as a discrete multi-criteria assessments (MCA) problem with \(n\) alternatives (firms) evaluated based on \(m\) minimization criteria (resources) and \(s\) maximization criteria (products). Each alternative has an observed score for each criterion, measured using a commensurable scale unit. The MCA method employs a mathematical algorithm to determine the weights of the criteria. The decision-maker ranks the alternatives according to their total weighted scores and selects the best one.

A typical MCA problem aims to evaluate each object decision-making unit (firm), denoted as \(DMU_o\), against \(n\) DMUs (firms) using \(m\) inputs (resources) and \(s\) outputs (products). Each EA approach employs a mathematical framework to assess the \textit{relative efficiency} of \(DMU_o\), pinpointing the necessary changes in inputs and outputs to attain a relative efficiency score of 1. Each DMU within the collection of DMUs is evaluated in turn as $DMU_o$.

While MCA problems rely on the observed dataset, their objectives differ. In this paper, the terms alternatives, minimization criteria, and maximization criteria are used interchangeably with DMUs, inputs, and outputs.

 \begin{table}[ht]
         \centering
     \caption{Numerical example of the MCA problems.} 
         \label{table:1}
     \centering
     \setlength{\tabcolsep}{2.5pt}
\renewcommand{\arraystretch}{0.95}
      \begin{tabular}{llcccccccc} 
      \hline\hline
\multicolumn{2}{c}{Criteria} & \multicolumn{6}{c}{Alternative-j (or $DMU_j$)} & \multicolumn{2}{c}{Decision variables}\\ 
\hline
\multicolumn{2}{c}{Minimization}& K & A & B & D & G & H & Adjustment& Virtual \\
\multicolumn{2}{c}{ vs. Maximization}&  &  &  &  &  &  & ratios & Unit prices\\ 
\hline
Input-$X_1$&$x_{1j}$(ton)   &1.6   &2.3   &1   &1.9   &1.8   &2.5  & $q_{1o}$& $v_{1o}(\$/ton)$ \\
Input-$X_2$ &$x_{2j}$(hr)   &145   &120   &29   &281   &250   &100 &$q_{2o}$&$v_{2o}(\$/hr)$  \\
\hline
 Output-$Y_1$&$y_{1j}(m^3)   $&1036   &1327   &567   &2446   &1794   &1000  &$p_{1o}$& $u_{1o}(\$/m^3$)  \\
 Output-$Y_2$ &$y_{2j}$($\%$)   &49   &97   &89   &97   &57   &70 &$p_{2o}$& $u_{2o}(\$/\%$) \\
\hline
\multicolumn{2}{c}{ } & \multicolumn{6}{c}{Decision variables} & \multicolumn{2}{c}{Systematic parameters}\\ 
\hline
\multicolumn{2}{c}{Variable Intensities}&$\pi_{oK}$&$\pi_{oA}$& $\pi_{oB}$ &
$\pi_{oD}$
&$\pi_{oG}$
&$\pi_{oH}$
& \multicolumn{2}{c}{$\kappa_o^1, \kappa_o^z, \kappa_o^2$, $\tau_o^\star(\$)$}\\
\hline\hline
 \end{tabular}
\end{table} 

\indent Since the 1960s, numerous multiple criteria decision-making (MCDM) methods have provided systematic frameworks for addressing complex decision-making problems across various fields \cite{1}. Sahoo and Goswami \cite{1} conducted a comprehensive review of MCDM methods, highlighting their advancements, applications, and potential future directions. Their study explored widely used MCDM methods, including AHP, TOPSIS, ELECTRE, PROMETHEE, DEA, MUAT, GRA, ANP, Fuzzy TOPSIS, and MOO. They noted that these methods are often constrained by subjective judgments and susceptible to individual biases. In our research, we argue that MCDM methods should limit the number of alternatives considered to ensure reliable and meaningful results. Decision-makers require a systematic approach to reduce the number of alternatives, mainly when subjective judgments play a critical role in the decision-making process.

\indent In specific scenarios, different MCDM methods yield significantly divergent rankings. Amor et al.\cite{2} analyzed 1,921 MCDM-related documents in English, published between 1990 and 2021. From this dataset, they selected 659 articles and 16 reviews as their final sample. They imported the corresponding BIB files containing complete records and cited references to perform bibliometric analysis using the bibliometrix tool. The statistical data were employed to map academic research on discrete multi-criteria decision analysis methods and their applications in sorting, classification, and clustering. Their study focused on nine widely used MCDM methods for practical applications. Regardless of the mapping approach, their findings revealed that existing MCDM methods often produce varying results, which can be a source of debate.

\indent Data Envelopment Analysis (DEA)\cite{3} \cite{4} is one of the major MCA methods. Introduced in 1978, DEA employs linear programming\cite{5} to estimate the adjustments and weights of the inputs and outputs of each decision-making unit (DMU), denoted as $DMU_o$. DEA's significant contributions have spurred extensive research in the field of MCA. Over the past 50 years, hundreds of textbooks and thousands of journal articles related to DEA have been published\cite{6}. Sickles and Zelenyuk\cite{7} provided a comprehensive review of MCA methods and the underlying assumptions of DEA theory. DEA models measure $DMU_o$'s relative efficiencies under the conditions of constant returns-to-scale (CRS) and variable returns-to-scale (VRS). The Additive Models have yielded incomplete solutions \cite{8}. The drawbacks of the original DEA models have been comprehensively discussed \cite{9} \cite{10} \cite{11} \cite{12} \cite{13} \cite{14}.

\indent The statistic-based Stochastic Frontier Analysis (SFA)\cite{15} is another foundational method in efficiency analysis (EA). SFA assumes that the estimation of the \textit{average production frontier} consists of two independent randomness components. One error term represents \textit{technical efficiency}, which is assumed to be negative due to $DMU_o$'s inability to reach the production frontier. The other error term accounts for unobservable and measurement difficulties. The deviation of $DMU_o$ from the mean efficiency can be calculated to evaluate the extent of its inefficiency. The average production frontier function is estimated by assuming that the errors follow a half-normal distribution and applying maximum likelihood estimation techniques. Madaleno and Moutinho\cite{16} discussed comparisons between DEA and SFA, noting that the reliability of the SFA model diminishes when the dataset does not satisfy its underlying assumptions. Over the past 50 years, thousands of articles and numerous textbooks have been published on SFA.

\subsection {Research Objectives} \label{sec:1.2}
\begin{figure}[ht] \centering
 \includegraphics
 [width=0.84\linewidth, height=4 in]{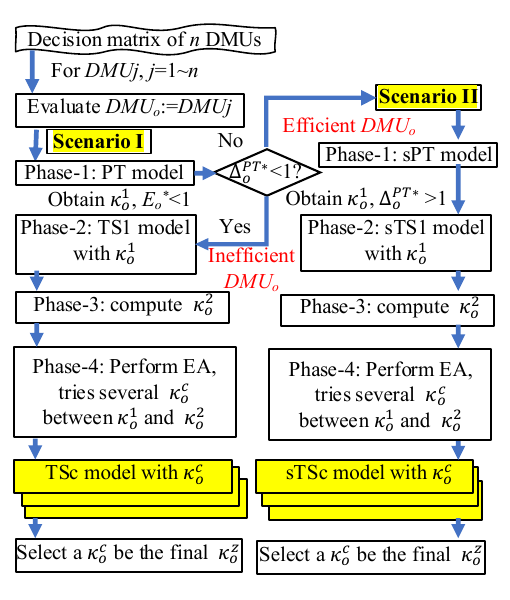} 
 \caption {The x b-VGA-EA method to assess DMUs.} \label{fig:1}
  \end{figure}

\begin{figure}[ht] 
\centering
 \includegraphics[width=0.84 \textwidth, height=4 in]{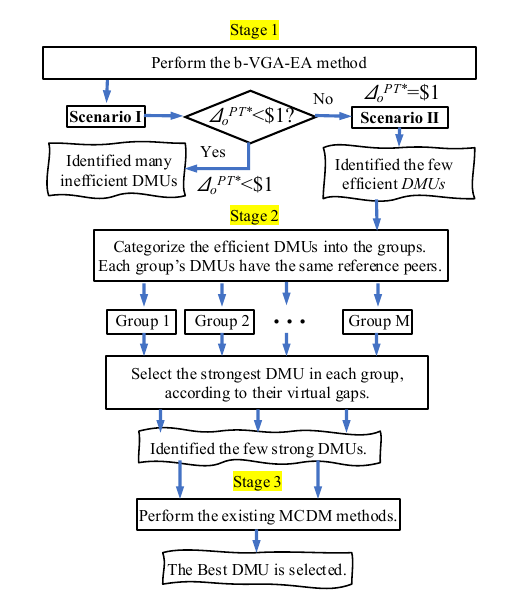} 
 \caption {The b-VGA-MCDM method to select the best alternative.} \label{fig:2}
  \end{figure}

This paper introduces a novel perspective and addresses the key limitations of traditional EA and MCDM methods. Figure \ref{fig:1} presents the flowchart of our new EA method, the \textit{Virtual Gap Analysis (VGA) procedure} based on best practices, abbreviated as b-VGA-EA. In this approach, each decision-making unit ($DMU_o$) evaluates itself against other DMUs to estimate the virtual gap.
\vspace{-1em}
Figure \ref{fig:2} illustrates our new MCA method, termed b-VGA, which integrates the b-VGA method to reduce the number of DMUs. This reduction allows for the efficient application of existing MCDM techniques to identify the best-performing DMU.

Finally, this paper explores EA and MCDM methods from the perspective of worst practices, where inputs and outputs are inverted. Using the MCDM method, the least efficient DMU can be effectively identified.

\subsubsection{Our EA method, b-VGA-EA method} \label{sec:1.2.1}
\indent In Figure \ref{fig:1}, Scenario I comprises four phases using the \textit{best-practice pure technical efficiency} (PT) and \textit{best-practice technical and scalar choice efficiency} (bTSc) VGA models. Scenario II includes the \textit{super-pure technical efficiency} (sPT) and \textit{super-technical and scalar choice efficiency} (sTSc) VGA models. The four linear VGA models are detailed in Sections \ref{sec:2} and \ref{sec:3}. A basic understanding of linear programming\cite{5} is required to comprehend the b-VGA-EA method fully.

\indent The four VGA models utilize decision variables (\textit{$q_{1o}$, $q_{2o}$, $p_{1o}$, $p_{2o}$}) and (\textit{$v_{1o}$, $v_{2o}$, $u_{1o}$, $u_{2o}$}) to represent the \textit{adjustment ratios} and \textit{virtual unit prices} of ($X_1$, $X_2$, $Y_1$, $Y_2$), respectively. The decision variables $\pi_{oK}$, $\pi_{oA}$, \dots, and $\pi_{oH}$ denote the \textit{intensities} of the six DMUs in evaluating \( DMU_o \).
The parameter \textit{unified goal price}, \( \tau_o^\star(\$) \), is systematically determined in two steps and serves as the common lower bound for the variable \textit{total virtual prices} of \( X_1 \), \( X_2 \), \( Y_1 \), and \( Y_2 \). The estimated optimal value of a decision variable in each VGA model is denoted with the superscript '$\star$.' The vector form of the estimated values are (\textit{$q_o^\star, p_o^\star$}), (\textit{$v_o^\star, u_o^\star$}), and $\pi_o^\star$. 

\indent In the bTSc and sTSc models, the additional sum of the intensities condition (SIC) requires the summation of variable intensities to equal a chosen scalar $\kappa_o^c$. The SIC corresponds to the dual variable $w_o$. Each model produces the systematic parameters \( \kappa_o^1 \) and \( \kappa_o^2 \). \( DMU_o \) selects \( \kappa_o^z \), a value between \( \kappa_o^1 \) and \( \kappa_o^2 \), to determine achievable adjustment ratios (\textit{$q_o^\star, p_o^\star$}) for optimizing its performance.

\indent In each model, a positive \( \pi_{oj}^* > 0 \) indicates that \( DMU_j \) is one of the \textit{reference peers} of \( DMU_o \). Each input (or output)  \( DMU_o \) is adjusted to the sum of the inputs (or outputs) of the reference peers, weighted by their respective intensities.

\indent Scenario I estimates the reduction of inputs ($q_o^\star$) and the expansion of outputs ($p_o^\star$). The PT model calculates the virtual gap of $DMU_o$, denoted as $\Delta_o^{PT\star}$, which represents the excess of \textit{virtual input} ($x_o v_o^\star$) over \textit{virtual output} ($y_o u_o^\star$). If $\Delta_o^{PT\star} < \$1$, $DMU_o$ is considered inefficient. In this scenario, $DMU_o$ proceeds to choose $\kappa_o^z$ in Phase-4 to achieve feasible values of $q_o^\star$ and $p_o^\star$. The resulting virtual gap is $\$0 \leq \Delta_o^ {bTSz\star} \leq \$1$, which is equal to the minimal value of $x_o v_o^\star$ - $y_o u_o^\star$ + $\kappa_o^ {bTSz} \times w_o^\star$. This reduction of the virtual gap to $\$0$ indicates that there is room for improvement in inputs and outputs by adjusting the ratios $q_o^\star$ and $p_o^\star$.

\indent Conversely, if $\Delta_o^{PT\star} = \$1$, $DMU_o$ is deemed efficient. In Scenario II, $DMU_o$ proceeds to choose $\kappa_o^z$ in Phase-4 to determine feasible $q_o^\star$ and $p_o^\star$ values. The virtual gap in this case is $\$1 \leq \Delta_o^{sTSz\star}$, which is equal to the minimal value of -$x_o v_o^\star$ + $y_o u_o^\star$ +$ \kappa_o^{sTSz} \times w_o^\star$. Reducing this virtual gap to $\$1$ reflects that $DMU_o$ possesses buffers for deterioration in inputs and outputs, which can be achieved by adjusting ratios $q_o^\star$ and $p_o^\star$ through expansion and reduction, respectively. 

\indent Table \ref{table:2} summarizes the solutions in assessing the inefficient \( DMU_K \). In column TS3-II, the optimal solutions of the decision variables for the TS3 VGA model can be found. \( DMU_B \) and \( DMU_D \) are the reference peers of \( DMU_K \), with intensities of 0.383 and 0.335, respectively.

\indent Figures \ref{fig:3} and \ref{fig:4} present the results of the assessments for the PT, TS1, and TS2 models, helping users better understand the b-VGA-EA method. The relationship between relative, technical, and scale efficiencies can be visualized in the 2D graphical representations. Similarly, the EA reports for the efficient \( DMU_B \) and \( DMU_D \) are shown on the right-hand side of Table \ref{table:2} and in the entire Table \ref{table:3}. 

\indent  Table \ref{table:3} summarizes the b-VGA-EA method reports in assessing the efficient \( DMU_B \) and \( DMU_D \). Figures \ref{fig:5} and \ref{fig:6} display 2D graphical representations of the assessments by the sTS1 and sTS2 models. These figures help visualize the relationship between relative, technical, and scale efficiencies. The technical and scale efficiencies have not been numerically computed at this stage.

\indent \( DMU_B \)'s relative efficiency score in the sTSz evaluation is 2.4479. According to the adjustment ratios  \( DMU_B \) in column sTSz-II, the ranking of strengths is as follows: \( Y_2 \), \( X_2 \), \( Y_1 \), and \( X_1 \), with \( X_1 \) being the weakest. These rankings differ from those observed in the sTS1 and sTS2 evaluations. The results of this numerical example demonstrate the depth of the EA method and its ability to provide insightful assessments.

 \subsubsection{Our MCDM method, b-VGA-MCDM method} \label{sec:1.2.2}
 \indent Figure \ref{fig:2} displays the b-VGA-MCDM method. In the first stage, the b-VGA-EA method is applied to identify the efficient DMUs, typically resulting in only a few efficient DMUs. These efficient DMUs, sharing the same reference peers, are categorized into groups in the second stage. Within each group, the efficient DMU with the highest virtual gap outperforms the others, and those with lower super-efficiency are excluded from consideration in the third stage. In the third stage, decision-makers can efficiently use existing MCDM methods to select the best DMU from the remaining candidates. They may also introduce several uncertain or qualitative criteria, as subjective judgment is often necessary. In some cases, when the number of candidates is minimal, decision-makers may not need to apply any existing MCDM methods.

\indent In our numerical example, $DMU_B$ and $DMU_D$
  have two different sets of reference peers. These DMUs are compared in the third stage. Decision-makers may assess the persistence capability of $DMU_B$ and $DMU_D$ with respect to their inputs and outputs. This study provides decision-makers with valuable insights and enables informed decision-making in dynamic selection environments.
  
\subsection{Early VGA Developments}\label{sec:1.3}
Hwang\cite{17} contributed to establishing CRS and VRS VGA models, where the former excludes and includes the SIC, equating to 1. However, issues arose concerning the virtual goal price  $DMU_o$ and the relationship between these models. The virtual goal price was subjectively assigned. Past research, including works by Liu and Liu\cite{18} \cite{19}, applied the CRS VGA model in two-phase production systems. 
The authors have provided consultation on unpublished theses that utilized these models. Specifically, Lin \cite{20} and Li \cite{21} applied the PT model in their respective thesis studies.

\subsection{Paper Organization} \label{sec:1.4}
This paper is organized as follows: Section 2 introduces the PT and bTSc VGA models for assessing efficiency. Section 3 discusses the sPT and sTSc VGA models for determining super-efficiency. Section 4 describes the four-phase b-VGA-EA method. Section 5 presents a numerical example to illustrate the two VGA scenarios. Finally, Section 6 concludes the paper and discusses the findings. Readers familiar with DEA may be interested in comparing DEA and VGA. To address this, we have included 26 footnotes to highlight how DEA models have inspired VGA models and to explain why DEA models cannot achieve certain aspects of efficiency analysis that VGA models can.

\section {Scenario-I: Evaluating DMUs' Inefficiencies} \label{sec:2} Scenario-I assesses the inefficiencies of DMUs by the PT and bTSc VGA models in a four-phase procedure.
\subsection {Notations} \label{sec:2.1}
 The following notations are used in VGA models. The set  $J$,  $I$, and  $R$, respectively, denotes the set of  $n$ DMUs,  $m$ inputs, and  $s$ outputs. The \textit {decision matrix}, denoted as  $(\textit{X}$, $\textit{Y})$, consists of column vectors $(x_j,y_j)$, in which the elements  $(x_{ij},y_{rj})$ represent positive \textit {observed volume} of each $Input_i$ ($X_i$) and $Output_r$ ($Y_r$) of each $DMU_j$.
 \footnote{\label{fn:1}     Cooper et al.\cite{4} introduced Equations (3.1) and (4.1) to represent the production possibility set (PPS) derived from the observed data of DMUs. The PPS without the convexity condition $e\lambda=1$ is referred to as the constant returns-to-scale (CRS) PPS, while the PPS with the convexity condition is known as the variable returns-to-scale (VRS) PPS. Each PPS is characterized by an efficiency frontier that envelops each $DMU_o$ at $(x_o, y_o)$.
 
    $P^{CRS}=\{(x_o,y_o)|x_o\leq X\lambda, y_o \geq Y\lambda, \lambda\geq 0.\}$
    
    $P^{VRS}=\{(x_o,y_o)|x_o\leq X\lambda, y_o \geq Y\lambda, (e\lambda=1,) \lambda\geq 0.\}$
    
\indent Cooper et al.\cite{4} presented the envelopment and multiplier programs of the Additive Model, given by Eqs. (4.34)–(4.38) and Eqs. (4.39)–(4.43), to estimate the vector of variable intensities $\lambda$.
 The model applies the production theory of neoclassical economics to evaluate \(DMU_o\) can be effectively characterized by the PPS. The PPS is defined by four key properties\cite{4}.

\indent Chapters 1 and 2 in Sickles and Zelenyuk\cite{7} provide twenty graphical illustrations of PPS in two dimensions, offering theoretical insights into CRS and VRS DEA models. However, a graphical representation of the efficiency frontier for PPS involving multiple inputs and outputs has yet to be developed. } 
 
  \indent  The decision variables ($q_{io}$, $p_{ro}$) 
  \footnote{\label{fn:2} Instead, DEA models used the slacks ($s_{io}^x, s_{ro}^y$) as decision variables, which have various measurement units.} 
  and ($v_{io}$, $u_{ro}$) 
  \footnote{\label{fn:3} DEA models defined ($v_{io}$, $u_{ro}$) as dimensionless weights.}  
  denote the dimensionless \textit{adjustment ratios} and \textit{virtual unit prices} (\$) of ($X_i$, $Y_r$), while their vectors are $(q_o, p_o)$ and $(v_o, u_o)$. $\pi_{jo}$ denotes the variable \textit{intensity} of $DMU_j$ in evaluating $DMU_o$, and ${\pi_o}$ is the vector form. 

  \indent We developed a \textit{two-step process} to determine a systematic \textit{unified goal price}  $\tau_o$, measured in the virtual currency \$. 
  \footnote{\label{fn:4} Instead, each DEA model developer defined the coefficient $b_i^x$ and $b_r^y$ be the the lower bound of weight $v_{io}$ and $u_{ro}$. The coefficients are called \textit{artificial goal weights}. }
  The superscripts  \# and   $\star$ indicate the optimal solutions for a decision variable in each VGA model, respectively, for Step I and Step II.
  
  \indent The symbol $w_o$ is the \textit{scalar unit price} (\$) corresponding to the SIC.
  \footnote{\label{fn:5} In the Additive Model, the symbol $\sigma_o$ is dimensionless and corresponds to the convexity condition, $ \sum_{\forall j\in J}\lambda_{oj}$=1. }
  $\kappa_o^c$ is the SIC scalar in the bTSc and sTSc VGA models.  $\omega_o^c$ is the \textit{virtual scalar price} (\$) of the SIC model $\kappa_o^c \times  w_o$. 
 
 \subsection{Pure Technical Efficiency (PT) model} \label{sec:2.2}

\textbf{The dual program} of the PT model. 
\begin{equation} \label{Eq:1}	
    \Delta_o^{PT\star} (\$)= \max_{q_o, p_o, \pi} \sum_{\forall i\in I} q_{io} \tau_o (\$) + \sum_{\forall r\in R}p_{ro} \tau_o (\$), (\forall o \in J);
 \end{equation}
\vspace{-1 em}
 \begin{equation} \label{Eq:2}
	s.t. \sum_{\forall j\in J} x_{ij} \pi_{jo} + q_{io} x_{io} =x_{io}  , \forall i\in I;
 \end{equation}
\vspace{- 1 em}
\begin{equation} 
 - \sum_{\forall j \in J} y_{rj} \pi_{jo} + p_{ro} y_{ro} = - y_{ro}, \forall r \in R;	
\label{Eq:3}
 \end{equation}
 \vspace{- 1.5 em}
\begin{equation} \label{Eq:4}
\pi_o, q_o  , p_o  \ge 0.	
 \end{equation}
Eq. (\ref{Eq:2}) and Eq. (\ref{Eq:3}) correspond to the primal program's variable unit prices, $v_{io}$ and $u_{ro}$.   Eq. (\ref{Eq:1}) measures the \textit{maximum total adjustment price} (TAP) of $DMU_o$. 
 
\indent \textbf{The primal program} of the PT model aims to measure the minimum \textit{total virtual gap} (TVG):
\footnote{\label{fn:6} Halická and Trnovská\cite{8} extensively reviewed DEA models and proposed a unified model to encompass various DEA approaches that address the Additive Model under VRS conditions, see footnote \ref{fn:1}. Without the four parentheses terms, it becomes the CRS Additive Model. \\
\textit{ The envelopment program}:

 \indent $F_o^\star= \max_{\pi_o\textbf, \textbf{s}_o^x, \textbf{s}_o^y} \sum_{\forall i\in I} s_{io}^x b_i^x  + \sum_{\forall r\in R} s_{ro}^y b_r^y, \forall o \in J;$	
	s.t. $\sum_{\forall j\in J } x_{ij} \lambda_{oj} + s_{io}^x   = x_{io}, \forall i \in I;$ 
	$- \sum_{\forall j\in J} y_{rj} \lambda_{oj} + s_{ro}^y  = - y_{ro}  , \forall r \in R;$ ($\sum_{\forall j \in J }\lambda_{oj}=1;$)
	$\lambda_{oj},\forall j \in J; s_{io}^x, \forall i \in I; s_{ro}^y, \forall r \in R  \ge 0.$ \\
   \textit{ The multiplier program}:
   
	$f_o^\star = \min_{\textbf{v}_o, \textbf{u}_o, \sigma_o}
	 \sum_{\forall i\in I} v_{io} x_{io}  - \sum_{\forall r\in R} u_{ro} y_{ro} (+ \sigma_o), \forall o \in J;
s.t. \sum_{\forall i\in I} v_{io} x_{ij} - \sum_{\forall r\in R} u_{ro} y_{rj} ( + \sigma_o ) \ge 0, \forall j \in J; 
v_{io} \ge b_i^x,\forall i \in I; 
u_{ro} \ge b_r^y, \forall r \in R;	
v_{io}, \forall i \in I; u_{ro}, \forall r \in R(, \sigma_o) free.$ }	

\footnote{\label{fn:7} 
See footnote \ref{fn:6}. The envelopment program of the DEA Additive Models evaluates dimensionless relative inefficiency as the total dimensionless slacks of inputs and outputs. However, a critical flaw lies in the fact that these slacks, $s_{io}^x$ and $s_{ro}^y$, actually possess distinct measurement units, contradicting the assumption of dimensional homogeneity. }
\footnote{\label{fn:8} Several DEA axioms fail to satisfy the fundamental conditions necessary for formulating linear programs. The decision variables, constraints, and objective functions of VGA models differ significantly from those of DEA Additive Models. Consequently, it is essential to abandon traditional DEA theory and instead formulate VGA models with both primal and dual programs to ensure their verification through duality principles. Unlike DEA models, VGA models define decision variables and coefficients in their actual measurement units, maintaining consistency with real-world data. }
 \indent The \textit{virtual gap condition} Eq. (\ref{Eq:6}) limits each $DMU_j$ to a total virtual gap, the excess of \textit{virtual input} and \textit{virtual output}, to a minimum of \$0. Eq. (\ref{Eq:7}) and Eq. (\ref{Eq:8}) restrict the \textit{virtual price} of $X_i$ and $Y_r$ to a unified lower bound \textit{unified goal price} $\tau_o$. 
 
 Eq. (\ref{Eq:6}), Eq. (\ref{Eq:7}), and Eq. (\ref{Eq:8}) correlate with the dual program's decision variable  $\pi_{jo}$, $q_{io}$, and $p_{ro}$, respectively. In Step I and Step II , substitute $\tau_o$  by  $\tau_o^\#$ = \(\$1\) and  $\tau_o^\star$   = $\$ \bar{t}$, respectively.

\begin{equation}
	\Delta_o^{PT\star} (\$)= \min_{v_o, u_o} \sum_{\forall i\in I} v_{io} x_{io}  - \sum_{\forall r\in R} u_{ro} y_{ro} , (\forall o \in J);	
 \label{Eq:5}
 \end{equation}
 \vspace{- 1.5 em}
 \begin{equation}
	s.t. \sum_{\forall i\in I} v_{io} x_{ij}  - \sum_{\forall r\in R}u_{ro}  y_{rj}  \ge 0(\$), \forall j\in J; \label{Eq:6}
 \end{equation}
 \vspace{- 1.5 em}
 \begin{equation}
	x_{io} v_{io} \ge\tau_o(\$), \forall i\in I;	
 \label{Eq:7}
 \end{equation}
 \vspace{- 1.5 em}
\begin{equation}\label{Eq:8}
	y_{ro} u_{ro} \ge \tau_o(\$), \forall r \in R;
 \end{equation}
 \vspace{- 1.5 em}
\begin{equation}
	v_o, u_o \quad free.	
 \label{Eq:9}
 \end{equation}
  \indent The optimal solutions of Eq. (\ref{Eq:1}) in Step I and Step II are depicted in Eq. (\ref{Eq:10}). Eq. (\ref{Eq:11}) expresses the multiple inputs and outputs of each $DMU_j$ in Step I and Step II are aggregated into the pairs of $(\alpha_j^\#, \beta_j^\#)$ and $(\alpha_j^\star, \beta_j^\star)$. In Step I, $\Delta_o^{PT\#}$ it may exceed that $\$1$.   
 \begin{equation}
 \begin{aligned}
&\$0 \leq \Delta_o^{PT\#} = \sum_{\forall i\in I} q_{io} ^\# \times \$1  + \sum_{\forall r\in R}p_{ro} ^\#  \times \$1 ; \\ 
&\$0 \leq \Delta_o^{PT\star } = \sum_{\forall i\in I} q_{io}^\star  \tau_o ^\star  + \sum_{\forall r\in R} p_{ro}^\star  \tau_o ^\star \leq \$1.
  \end{aligned}  \label{Eq:10} 
 \end{equation} 
 \begin{equation} \label {Eq:11}
\begin{aligned}
\$0 \leq (TVG)=(vInput)-(vOutput) ;\\
\$0 \leq  \Delta_j^ {PT\#}  = v_o^\# x_j - u_o^\# y_j = \alpha_j^\# - \beta_j^\# ; \\
\$0 \leq  \Delta_j^ {PT\star }  = v_o^\star  x_j - u_o^\star  y_j =\alpha_j^\star  - \beta_j^\star \leq \$1.
\end{aligned}
 \end{equation}

\indent  \textbf{DMU-o determines its unified goal price $\tau_o^\star$} to ensure $ \Delta_o^{PT\star}$ and $\Delta_o^{PT\star}$ in Step II are within the range of ($\$0, \$1$). The relationship
$\tau_o^\# :\tau_o^\star$   =
$ 1 :  \bar{t}$ =  $\Delta_o^{PT\#}:\Delta_o^{PT\star}$ could be derived from the \textit{relation equation} $\bar{t}  (   \alpha_o^\#- \beta_o^\# ) = 1(\alpha_o^\star-\beta_o^\star$).  Use
 Eq. (\ref{Eq:12}) to obtain the dimensionless value of  $\bar{t}$. 
\begin{equation}
	\$\bar{t} = \$1/\alpha_o^\# \textrm {  and  } \tau_o^\star   = \$ \bar{t}. 
	\label{Eq:12}
  \end{equation}
  Dividing the relation equation by $\alpha_o^\star$ and substituting $\bar {t}$, it becomes $ (\$1/\alpha_o^\star)( 1- \beta_o^\#/\alpha_o^\# ) =(1-\beta_o^\star/\alpha_o^\star )$. 
  
  \indent Because of Eq. (\ref{Eq:5}) obtained $ \$0 \leq( \alpha_o^\#- \beta_o^\# )$ and $ \$0\leq (\alpha_o^\star-\beta_o^\star) $, we have $0\leq( 1- \beta_o^\#/\alpha_o^\# )\leq 1$ and $0\leq (1-\beta_o^\star/\alpha_o^\star )\leq 1$.   Using $\tau_o^\star   = \$ \bar{t}$ in Step II, one should have the solutions  $\alpha_o^\star =\$1$ and $\beta_o^\star \leq \$1$ to confirm the relation equation within the range of (0, 1).  $\Delta_o^{PT\star}$ should be within the range of ($\$0, \$1$).
  
  \indent Therefore, the relative efficiency $E_o^ {PT\star}$ is ensured between 0 and 1, as shown in Eq. (\ref{Eq:13}).

\begin{equation}
\begin{aligned}\label{Eq:13} 
	0< E_o^ {PT}  = \beta_o^\star/\alpha_o^\star \le 1.  
\end{aligned}
\end{equation}

\indent Certainly, normalizing Step I solutions to attain Step II  solutions is a common practice in mathematical optimization problems.  Eq. (\ref{Eq:14}) likely represents this normalization process, where the Step I solutions are adjusted or transformed to achieve the solutions in Step II. This normalization might involve scaling or modifying Step I solutions to fulfill specific conditions required for Step II.
\begin{equation} 
  (\Delta_o^{PT\star },  \Delta_o^{PT\star },  v_o^\star, u_o^\star) = \bar{t} \times (\Delta_o^{PT\#}, \Delta_o^{PT\#}, v_o^\#, u_o^\#).
\label {Eq:14}
\end{equation}

\indent The symbol  $\mathcal{E}_o^ {PT}$  denotes the set of \textit{reference peers}  $DMU_o$ in the PT evaluations. $DMU_j$ belongs to $\mathcal{E}_o^ {PT}$ has $\pi_{jo}^\star>0$ and is an efficient one.  At the same time, the estimated intensity of the other inefficient $DMU_j$, $\pi_{jo}^\star =0$. $DMU_o$ estimates the benchmark of each $X_i$ and $Y_r$,  $\widehat{x}_{io}^ {PT\star}$ and $\widehat{y}_{ro}^ {PT\star}$ via Eq. (\ref{Eq:15}). $DMU_o$ mimics reference peers with their estimated intensities, $\pi_{jo}^\star$, $\forall j\in \mathcal{E}_o^ {PT}$.

\footnote {\label{fn:9} Each DEA model used the same set of artificial goal weights to identify the reference peers for every $DMU_o$. The set of reference peers and $DMU_o$ could have heterogeneous configurations of the observed dataset for more details \cite{20} \cite{23} \cite{24}. To solve a VGA model, $DMU_o$ systematically determines its unified goal price $\tau_o^\star$ in two steps. It always identifies its homogeneous reference peers. }
\footnote{\label{fn:10} Halická and Trnovská\cite{8} listed objective functions and artificial goal weights used in several DEA Additive Models in  Tables 2 and  3. The Additive CRS Models may obtain infeasible relative inefficiency $F_o^\star$, which is not within the range (0,1) because the estimated $x_o^\star v_o^\star \neq{1} $. However, the artificial goal weights in each model could not ensure that $x_o^\star v_o^\star =1.$ The Additive Models are unreliable and cannot be used. }
 \footnote {\label{fn:11} The slack-based-measure (SBM) model is a special one among the Additive Models; SBM models always have satisfactory results $0 \leq F_o^\star (=f_o^\star) < 1 $. Therefore, SBM is the most frequently used method for EA.
 
\indent Using a numerical example with multiple inputs and outputs to test those SBM models would find the result $F_o^\star =f_o^\star
 >(1- y_o u_o^\star / x_o v_o^\star)$, the relative efficiency score $ y_o u_o^\star/ x_o v_o^\star$ was underestimated. Only if  $x_o v_o^\star=1$ could satisfy the condition $F_o^\star =f_o^\star
 =(1- y_o u_o^\star / x_o v_o^\star)$. 
 
 \indent In SBM literature, the solutions $v_o^\star$ and $u_o^\star$ usually are not shown. The inaccurate SBM estimations will not be realized.
 For the original SBM models \cite{25} \cite{26}, the development program aims to measure the efficiency score, which is a fractional expression. That led to the artificial goal weights in the multiplier program, $v_{io}\geq b_i^x [= 1/(mx_{io})]$  and $ u_{ro} \geq b_r^y [=(x_o v_o - y_o u_o)/ (s y_{ro})]$. Apparently, SBM models have inaccurate evaluations. }
\begin{equation}
\begin{aligned}
	&\widehat{x}_{io}^ {PT\star} = \sum_{\forall j\in \mathcal{E}_o^ {PT}} x_{ij}  \pi_{jo}^\star  = x_{io} (1 - q_{io}^\star ), \forall i\in I; \\
	&\widehat{y}_{ro} ^ {PT\star} = \sum_{\forall j\in \mathcal{E}_o^ {PT}}y_{rj} ,\pi_{jo}^\star =y_{ro} (1 + p_{ro}^\star  ), 
 \forall r\in R.
\end{aligned}
\label {Eq:15}
\end{equation}
Assessing $DMU_o$ with the adjusted $X_i$ and $Y_r$, $\widehat{x}_{io}^ {PT\star}$ and $\widehat{y}_{ro} ^ {PT\star}$, as shown in  Eq. (\ref{Eq:15}), will have the PT efficiency $\widehat{E}_o^ {PT}$ equals 1. Let the total of estimated intensities be the first SIC scalar, $\kappa_o^1$, computed via Eq. (\ref{Eq:16}). 
\begin{equation} \label{Eq:16}
    \sum_{\forall j\in \mathcal{E}_o^ {PT}} \pi_{jo}^\star= \kappa_o^1. 
  \end{equation}
 
 \indent \textbf{Numerical example. }In the columns of PT-I and PT-II of Table \ref{table:2}, summarized the solutions of Step I and Step II.  The optimal values of $\tau_o$,  $(v_o, u_o)$, $(v_{io}, u_{ro})$, and $(\alpha_j, \beta_j)$ are altered, while the other are unchanged. $DMU_B$ and $DMU_D$ are belong to $\mathcal{E}_o^ {PT}$, while their intensities $\pi_{oB}^\star$ and $\pi_{oD}^\star$ are positive. The first SIC scalar is computed.

\subsection {Technical and Scalar Choice (TSc) Model} \label{sec:2.3}
Adding Eq. (\ref{Eq:20}) to the PT-TAP program will have the bTSc-TAP program.

\footnote{\label{fn:12} In DEA VRS models, the \textit{convexity condition} is expressed as \(\sum_{\forall j \in J}\lambda_{jo} = 1\) 
 (see footnote \ref{fn:1}). While the CRS PPS polyhedral space does not include this convexity condition, the VRS PPS space incorporates it. The axioms of DEA theories are well established\cite{4} \cite{7}. Based on these theories and their graphical representations, each observed point corresponding to a \(DMU_o\) can adjust its single input and single output either by reduction or expansion. The resulting slacks, 
\(s_{1o}^x\) and \(s_{1o}^y\) project the point onto the efficient frontier of the envelope.

\indent The envelopment analysis has been extended to the PPS, where DMUs operate with multiple inputs and outputs. This extension constitutes one of the fatal theoretical errors of DEA. The envelopment program of the Additive Model (see footnote \ref{fn:6}) employs decision variables of slacks
 \(s_{io}^x\) and \(s_{ro}^y\) to measure the maximum relative inefficiency of \(DMU_o\). \indent Based upon the DEA theory, the decision variables of intensities of DMUs $\lambda_{oj}, \forall j \in J$ in analyzing the PPS space are the same in analyzing the Additive Models\cite{4} \cite{7}. However, the two sets of intensities should be different. The Additive Model constructs a convex, feasible polyhedral space using hyperplanes and program conditions. In VGA models, use $\pi_o$ as the variable intensities of DMUs.
\indent In the CRS (VRS) envelopment program, the absence (presence) of the condition \(\sum_{\forall j \in J}\lambda_{jo} \)=1 distinguishes the CRS space from the VRS space. Importantly, this condition is not a convexity constraint within the convex space but rather a simple additional constraint. Adding this constraint reduces the CRS space to the VRS space. Implying the PPS convexity condition to the envelopment program introduces another fatal theoretical error in DEA.

\indent The spaces defined by the PPS and the envelopment program differ significantly. The polyhedral PPS space comprises a finite set of points corresponding to the observed data. The reference peers of \(DMU_o\) are points situated on the piecewise linear PPS frontier, which cannot be represented either graphically or mathematically as a production function. In contrast, the feasible space of the envelopment program is a convex polyhedral space constructed using hyperplanes.

\indent The optimal solution of the envelopment program is located at an extreme point, formed by the intersection of multiple hyperplanes and the objective function. Conversely, the projection point in the PPS space lies on the frontier, where the efficiency reference peers are situated. These points are distinct and are determined based on different criteria.

\indent The envelopment programs for CRS Additive Models yield incomplete evaluations (see footnotes \ref{fn:10} and \ref{fn:11}), as do those for VRS Additive Models \cite{8}.}
Eq. (\ref{Eq:20})  corresponds to the free-in-sign decision variable, $ w_o (\$)$. Eqs. (\ref{Eq:18}),(\ref{Eq:19}), \ref{Eq:20}, and (\ref{Eq:21}) construct the feasible region of the linear programming. 

\noindent \textbf {bTSc model, TAP (dual) program:}
\begin{equation}\label {Eq:17} 
	\Delta_o^{bTSc\star}(\$)= \max_{q_o, p_o} \sum_{\forall i\in I} q_{io} \tau_o + \sum_{\forall r\in R}p_{ro}   \tau_o, (\forall o \in J);	
\end{equation}
 \vspace{- 1 em}
 \begin{equation} \label{Eq:18}
	s.t. \sum_{\forall j\in J } x_{ij} \pi_{jo} + q_{io} x_{io}  = x_{io}, \forall i \in I;	
 \end{equation}
  \vspace{-1 em}
\begin{equation}\label {Eq:19}
	- \sum_{\forall j\in J} y_{rj} \pi_{jo} + p_{ro} y_{ro} = - y_{ro}  , \forall r \in R;	
\end{equation}
 \vspace{-1 em}
 \begin{equation}\label {Eq:20}
	 \sum_{\forall j\in J } \pi_{jo}  = \kappa_o^c;
  \end{equation}
   \vspace{-1.5em}
  \begin{equation}\label {Eq:21}
	\pi_o,  q_o, p_o  \ge 0.	
\end{equation}
  \noindent \textbf {bTSc model, TVG (primal) program:}

\begin{equation}\begin{aligned}
	\Delta_o^{bTSc\star} (\$)=  \min_{x_o, v_o,  w_o }  \sum_{\forall i\in I} v_{io} x_{io}  &- \sum_{\forall r\in R} u_{ro} y_{ro} + \kappa_o^c  w_o, (\forall o \in J);
 \label {Eq:22}
 \end{aligned}
 \end{equation}
  \vspace{-1.5 em}
 \begin{equation}\label {Eq:23}
s.t. \sum_{\forall i\in I} v_{io} x_{ij} - \sum_{\forall r\in R} u_{ro} y_{rj}  + 1  w_o  \ge 0(\$), \forall j \in J; 
\end{equation} 
 \vspace{-1.5 em}
 \begin{equation}\label{Eq:24}
x_{io} v_{io} \ge\tau_o(\$),\forall i \in I; 
\end{equation}
 \vspace{-1.5em}
 \begin{equation}\label {Eq:25}
y_{ro} u_{ro} \ge \tau_o(\$), \forall r \in R;	
 \end{equation}
  \vspace{-1.5em}
 \begin{equation}\label {Eq:26}
v_o, u_o,  w_o \quad free.
\end{equation}

\indent In  Eq. (\ref{Eq:22}), the elements $(\kappa_o^c  w_o^\#)$ and \((\kappa_o^c  w_o^\star) \) are repressed by the symbols  \(\omega_o^{c\#}\) 
and  \(\omega_o^{c\star}\), the \textit{ virtual scalar \$  (vScalar)} in Steps I and II.
 
As depicted in  Eq. (\ref{Eq:27}),  \(\gamma_o\) and \((1- \gamma_o)\) are the proportions of the total adjustment prices of inputs and outputs.
\footnote{\label{fn:13} The multiplier program in each of the DEA Additive VRS Models, as noted in footnote \ref{fn:6}, defines the optimal dimensionless relative inefficiency 
$f_o^*$ as comprising three components: $x_o v_o^\star -y_o u_o^\star + 1 \times w_o^\star$. However, the role of $1\times w_o^\star$ is not clearly defined. The convexity condition of the PPS is applied as the convexity condition for the feasible space of the Additive Model. This model assumes that 
$f_o^\star$ represents inefficiency. Nevertheless, the model encounters infeasible solutions\cite{8}, as discussed in footnotes \ref{fn:10} and \ref{fn:11}.

\indent Instead, the minimized TVG $\Delta_o^{bTSc\star}$ in Eq. (\ref{Eq:23}) is the virtual gap ($\$$). We partitioned TVG into two parts. }
If \(w_o^\#\) (or \( w_o^{\star}\)) equals \(\$ 0\), let $\gamma_o$ equal 0.5. Partition the \textit{vScalar} into two parts to reflect the effects of the SIC on inputs and outputs of $DMU_o$. We denote them as \emph{vScalar in inputs (ivScalar)} and \emph{vScalar in outputs (ovScalar)}; \textit{vScalar} equals \textit{ivScalar} plus \textit{ovScalar}.
\begin{equation}\begin{aligned}
	\gamma_o : (1-\gamma_o) = \tau_o^\# \times\sum_{\forall i\in I}q_{io}^\# :\tau_o^\# \times\sum_{\forall r\in R}p_{ro}^\# = \tau_o^\star  \times\sum_{\forall i\in I}q_{io}^\star  :\tau_o^\star \times\sum_{\forall r\in R}p_{ro}^\star .
 \end{aligned}\label {Eq:27}
 \end{equation}
\indent Eq. (\ref{Eq:28}) is the expression of Eq. (\ref{Eq:22}) in Step II, and the definitions of $\alpha_o^{c\star}$ and $ \beta_o^{c\star}$. The symbols $\alpha_o^{c\#}$ and $ \beta_o^{c\#}$ represent the two parts of $\Delta_o^{bTSc\#}$ in Step I. Usually, $\Delta_o^{bTSc\#}$ may exceed $\$1$.  We gave several names for the mathematical terms in the equation.  
 \begin{equation}\begin{aligned}\label {Eq:28}
\$0 \leq \Delta_o^{bTSc\star} &= v_o^\star x_o-u_o^\star y_o+\omega_o^{c\star } = [v_o^\star x_o+(1-\gamma_o)\omega_o^{c\star }]-(u_o^\star   y_o-\gamma_o\omega_o^{c\star }) \\
&=[gvInput+ivScalar]-[gvOutput-ovScalar]\\
&= avInput-avOutput  = \alpha_o^{c\star} - \beta_o^{c\star}  \leq \$1; \\
\$0 \leq \Delta_o^{bTSc\#} &= v_o^\# x_o-u_o^\# y_o+\omega_o^{c\#} = [v_o^\# x_o+(1-\gamma_o)\omega_o^{c\# }]-(u_o^\#   y_o-\gamma_o\omega_o^{c\# }) \\
&= avInput-avOutput =\alpha_o^{c\#} - \beta_o^{c\#}.
 \end{aligned} \end{equation}

 \indent Similarly, Eq. (\ref{Eq:29}) expressed the solutions of Eq. (\ref{Eq:23}) in Step II are converted into two parts for each $DMU_j$, $\alpha_j^{c\star }$ and $\beta_j^{c\star }$.  The estimated \textit{vScalar}, \(1w_o^\star\) is decomposed into two components. The exact conversions are applied to the solutions of Step I. 
  \begin{equation}\begin{aligned} \label{Eq:29}
\$0 \leq \Delta_j^ {bTSc\star} &= v_o^\star x_j-u_o^\star y_j+ 1 w_o^{\star}\\
&= [v_o^\star x_j+(1-\gamma_o) 1 w_o^{\star}]-(u_o^\star   y_j-\gamma_o 1 w_o^{\star})
= \alpha_j^{c\star } - \beta_j^{c\star }, \forall j \in J.
 \end{aligned} \end{equation}
 Each $DMU_j$ is expressed by the pair of \emph{(avInput, avOutput)}, the virtual scales.

\indent  \textbf{DMU-o determines its unified goal price, $\tau_o^\star$}, to ensure $ \Delta_o^{bTSc\star}$ and $\Delta_o^{bTSc\star}$ are within the range of (\$0 and \$1). Eq. (\ref{Eq:28}) supplied ($\alpha_o^{c\#}, \beta_o^{c\#}$) and ($\alpha_o^{c\star }, \beta_o^{c\star}$).  The relationship 
$\tau_o^\# :\tau_o^\star  =
1 :  \bar{t} = \Delta_o^{bTSc\#}:\Delta_o^{bTSc\star} $ 
could be derived from the \textbf{relation equation} $ (  \bar{t}  \alpha_o^{c\#}- \bar{t} \beta_o^{c\#} ) = (\alpha_o^{c\star}-\beta_o^{c\star}$).  Use
 Eq. (\ref{Eq:12}) to obtain the dimensionless value of  $\bar{t}$. 
\begin{equation} \begin{aligned} \label {Eq:30}
	\bar{t }= \$ 1 / \alpha_o^{c\#}=  \$1 / [v_o^\#  x_o +(1 - \gamma_o)\omega_o^{c\#}] .
\end{aligned} \end{equation}
 Dividing the relation equation by $\alpha_o^{c\star}$ and substituting $\bar{t}$, we have $ (  \$1/\alpha_o^{c\star})(1- \beta_o^{c\#}/ \alpha_o^{c\#} ) = (1-\beta_o^{c\star}/\alpha_o^{c\star}$).  
 
 \indent Because of Eq. (\ref{Eq:28}) obtained $ \$0 \leq( \alpha_o^{c\#}- \beta_o^{c\#} )$ and $ \$0\leq (\alpha_o^{c\star}-\beta_o^{c\star}$), we have $0\leq( 1- \beta_o^{c\#}/\alpha_o^{c\#} )= (1-\beta_o^{c\star}/\alpha_o^{c\star} ) \leq 1$.   Using $\tau_o^\star   = \$ \bar{t}$ in Step II should have the solutions $\alpha_o^{c\star} =\$1$ and $\beta_o^{c\star} \leq \$1$ to confirm the relation equation within the range of (0,1). Then, $DMU_o$ could ensure $\$0 \leq \Delta_o^{bTSc\star}\leq \$1$. $DMU_o$ systematically determines the virtual unit prices ($v_o^\star, u_o^\star$) and virtual scales ($\alpha_j^{c\star}, \beta_j^{c\star}$) for each $DMU_j$. 
 \footnote{\label{fn:14} In the Additive CRS and VRS Models, the estimated weights ($v_o^\star, u_o^\star$) are biased by the artificial goal weights  ($b_i^x, b_r^y$). In DEA literature, the virtual scales ($\alpha_o^{c\star}, \beta_o^{c\star}$) were not applied. However, the solutions ($v_o^\star, u_o^\star$) are not revealed.} 
 Therefore, as shown in Eq. (\ref{Eq:31}), the relative efficiency $E_o^ {bTSc\star}$ is ensured between 0 and 1.
 \footnote{\label{fn:15} The DEA Additive VRS Models have infeasible solutions \cite{8}, the reasons are addressed in footnotes \ref{fn:10}, \ref{fn:11}, and \ref{fn:13}. Employing the PT model could obtain comprehensive solutions to the DEA CRS model. Using $\kappa_o^ {bTSc}=1$ to the bTSc model should have the DEA VRS comprehensive solutions.}
\begin{equation}
\begin{aligned}
 0 < E_o^ {bTSc} &=\frac{
	 u_o^\star y_o-\gamma_o \omega_o^{c\star}}{v_{o}^\star x_o+(1-\gamma_o)\omega_o^{c\star} } = \beta_o^{c\star}/\alpha_o^{c\star} \leq 1.	\end{aligned}
 \label {Eq:31}
  \end{equation}

\indent Normalizing the solutions of Step I by the dimensionless value $\bar{t}$ would obtain the solutions of Step II, as shown in  Eq. (\ref{Eq:32}). 
\begin{equation}\begin{aligned}\label {Eq:32}
	(\Delta_o^{bTSc\star }, \Delta_o^{bTSc\star }, v_o^\star , u_o^\star ,  w_o^{\star})
	= \quad\bar{t}\times(\Delta_o^{bTSc\#},\Delta_o^{bTSc\#} ,v_o^\#, u_o^\# , w_o^{\#}).
\end{aligned}
\end{equation}
\indent Use  Eq. (\ref{Eq:33}) to compute the benchmark of each $X_i$ and $Y_r$,  \(\widehat{x}_{io}^ {bTSc\star}\) , and  \(\widehat{y}_{ro}^ {bTSc\star}\). $DMU_o$ imitates the reference peers with their estimated intensities, $\pi_{jo}^\star $. $DMU_o$ would become efficient with the benchmarks. 
\begin{equation}
\begin{aligned}\label {Eq:33} 
	\widehat{x}_{io}^ {bTSc\star} &= \sum_{\forall j\in \mathcal{E}_o^ {bTSc}}x_{ij}\pi_{jo}^\star = x_{io}(1 - q_{io}^\star ), \forall i\in I;\\
	\widehat{y}_{ro} ^ {bTSc\star} &= \sum_{\forall j\in \mathcal{E}_o^ {bTSc}}y_{rj}\pi_{jo}^\star = y_{ro}(1 + p_{ro} ^\star ), \forall r\in R .   \end{aligned} \end{equation}

\indent Assess $DMU_o$ with \(\widehat{x}_{io}^ {bTSc\star}\) and \(\widehat{y}_{ro} ^ {bTSc\star} \) by the bTSc model will have $\Delta_o^{bTSc\star}=\$0$, and the bTSc efficiency   $\widehat{E}_o^ {bTSc\star }$ equals 1.

 \indent The symbol  \(\mathcal{E}_o^ {bTSc}\) denotes the set of reference peers in evaluating $DMU_o$ by the bTSc model, where each reference $DMU_j$ has \(\pi_{jo}^\star >0\). The other inefficient $DMU_j$ has \(\pi_{jo}^\star =0\). In  Eq. (\ref{Eq:29}), any best $DMU_j$ belongs to \(\mathcal{E}_o^ {bTSc}\) has  $\alpha_j^{c\star } = \beta_j^{c\star }$, while the other remaining DMUs have  $\alpha_j^{c\star } > \beta_j^{c\star }$. 

\indent \textbf{Numerical example. }In the columns of TS1-I and TS1-II of Table \ref{table:2}, summarize the solutions in Step I and Step II. The first SIC scalar obtained in the PT model is used. The optimal values of $\tau_o$  $(v_o, u_o)$, $(v_{io}, u_{ro})$, $w_o$, $\omega_o^c$ and $(\alpha_j, \beta_j)$ are altered, while the others are unchanged. The homogeneous reference peers $DMU_B$ and $DMU_D$ belong to $\mathcal{E}_o^{TS1}$, while their intensities $\pi_{oB}^\star$ and $\pi_{oD}^\star$ are positive. The comparisons between the columns of PT-II and TS-II are necessary.
\footnote{\label{fn:16} \indent SFA methods may suffer from the heterogeneous configurations of the dataset, making it a challenging task to identify outliers. SFA estimates the average production frontier for each dependent (output, $Y_r$) and \textit{m} independent (inputs) of the \textit{n} DMUs. In contrast, VGA models assess 
$DMU_o$ with \textit{m} inputs and \textit{s} outputs simultaneously, considering its reference peers rather than relying on the average production frontier. }

\indent In the columns of TS2-I and TS2-II of Table \ref{table:2} , $DMU_o$ used the second SIC scalar and obtained the solutions in Phase-3. Comparing with the columns of TS1-I and TS1-II, we should find the solutions that have changed.

\subsection{Duality Properties}\label{sec:2.4}
We formulated VGA models based on the theories of linear programming\cite {4}. Each money coin has two faces, while every linear programming model has primal and dual programs. Dualities are used to examine the correctness of the model formulations. In particular, the measurement units used for the parameters and decision variables should be consistent in the primal and dual programs.
\footnote{\label{fn:17}In DEA literature, many proposed DEA models do not offer the multiplier program and its solutions. As a result, the dualities of those models cannot be examined. Additionally, the formulations of these DEA models are dimensionless \cite{27}. } 
When the optimal TVG and TAP of the PT and bTSc models are equal, as stated in Eq. (\ref{Eq:34}), it suggests an equilibrium or balance between the defined conditions or objectives within these models. This equality could signify an alignment between the total virtual gap and total slack price considerations, indicating a harmonized solution or relationship between the gap and price metrics.
\begin{equation}\begin{aligned} \label{Eq:34}
\Delta_o^{PT\star}(\$)=\Delta_o^{PT\star} (\$)\quad \textrm {and}   \quad
\Delta_o^{bTSc\star}(\$)=\Delta_o^{bTSc\star}(\$).
  \end{aligned} \end{equation}
\indent The strong complementary slackness conditions of the two VGA models are shown in Eq. (\ref{Eq:35}) and Eq. (\ref{Eq:36}).
\begin{equation} \begin{aligned} \label {Eq:35}
\relax [\sum_{\forall j\in \mathcal{E}_o^{VGA}} x_{ij} \pi_{jo}^\star - x_{io} (1 - q_{io}^\star)] v_{io}^\star =\$0,
\forall i\in I.  
\end{aligned}  \end{equation}
 \vspace{-1em}
\begin{equation} \begin{aligned}
\relax [\sum_{\forall j\in \mathcal{E}_o^{VGA}} (y_{rj}\pi_{jo}^\star )- y_{ro}(1+p_{ro}^\star )]u_{ro}^\star =\$0,
 \forall r\in R .
 \label {Eq:36} \end{aligned} \end{equation}
\indent $v_o$ and $u_o$ are free-in-signs in the TVG programs. Since \((\textit{X, Y})>0\) and \(\tau_o^\star >\$0\) yields the estimated virtual unit prices, \(v_o^\star  >0, u_o^\star >0\). Therefore, the braces at the left-hand side of  Eq. (\ref{Eq:35}) and Eq. (\ref{Eq:36}) equal zero. In other words, $DMU_o$ does not need to perform the sensitivity analysis on the parameters in the model ($x_o, y_o$). The VGA models obtained robust estimations of $p_{io}^\star$ and $q_{ro}\star$ so that $DMU_o$ ensued to project on the \textit{b-Efficiency Equator}, see Figures \ref{fig:3} through \ref{fig:6}.
\footnote{\label{fn:18} Using each existing DEA Additive CRS and VRS Model to evaluate each $DMU_o$ in the PPS with multiple inputs and outputs should provide a similar 2D graphical representation to demonstrate the solutions. The adjusted $DMU_o$ should have a relative efficiency score of 1 and should project onto the b-Efficiency Equator. However, the 2D figure may reveal that the DEA model provides incomplete evaluations. In contrast, VGA models have been verified to offer comprehensive evaluations, as demonstrated in Figures \ref{fig:3} through \ref{fig:6}.}

\indent For the PT model, each $DMU_j$ has the property as shown in  Eq. (\ref{Eq:37}). When $(v_o^\star x_j - u_o^\star y_j )=0$, $DMU_j$ belongs to \(\mathcal{E}_o^ {PT}\)  and \(\pi_{jo}^\star >0\). The \emph{b-Efficiency Equator} is the diagonal line at the origin in Figure \ref{fig:3}, \((v_o^\star x_j - u_o^\star y_j)=\$0, \forall j\in \mathcal{E}_o^ {PT}\).
\begin{equation}
	(v_o^\star x_j - u_o^\star y_j )\pi_{jo}^\star = \$0, \forall j\in J .
	\label{Eq:37}
\end{equation} 
For the bTSc model, each $DMU_j$ has the property shown in  Eq. (\ref{Eq:38}). $w_o^\star$ which is the estimated unit price of the SIC scalar and $\kappa_o^c w_o^\star$  is the \textit{vScalar}. When $(v_o^\star x_j - u_o^\star y_j + 1 w_o^\star ) = \$0$, $DMU_j$ belongs to $\mathcal{E}_o^ {bTSc}$ and $\pi_{jo}^\star>0$. Otherwise, $(v_o^\star x_j - u_o^\star y_j + 1  w_o^\star ) > \$0$, $DMU_j$ does not belong to $\mathcal{E}_o^ {bTSc}$ and $\pi_{jo}^\star=0$. The b-Efficiency Equator is the diagonal line at the origin in Figure \ref{fig:4}, $(v_o^\star x_j - u_o^\star y_j + 1  w_o^\star)=\$0, \forall j\in \mathcal{E}_o^ {bTSc}$. 
 \vspace{-.5em}
\begin{equation}
	 (v_o^\star x_j - u_o^\star y_j + w_o^\star ) \pi_{jo}^\star = \$0, \forall j\in J.
	 	\label{Eq:38} \end{equation}
 Eq. (\ref{Eq:39}) depicts the relationships between the estimated decision variables in the TVG and TAP programs.
\begin{equation}\begin{aligned}
	(v_{io}^\star x_{io} - \tau_o^\star) q_{io}^\star = \$0,\forall i \in I ; (u_{ro}^\star y_{ro} - \tau_o^\star) p_{ro}^\star= \$0,\forall r\in R .\end{aligned} \label{Eq:39} \end{equation}
 In  Eq. (\ref{Eq:40}), $\sum_{\forall j\in \mathcal{E}_o^ {bTSc}}\pi_{jo}^\star - \kappa_o^c = 0$ the TVG of the bTSc model is decreased, constant, and increased when $w_o^\star>\$0,  w_o^\star=\$0, \textrm{ and } w_o^\star<\$0$respectively. 
\begin{equation} \label{Eq:40}
	(\sum_{\forall j\in \mathcal{E}_o^ {bTSc}} \pi_{jo}^\star -\kappa_o^c ) w_o^\star  = \$0.			\end{equation}
    
\subsection{Post-Analysis of the Optimal Solutions}\label{sec:2.5}
The optimal solutions from the PT and bTSc models have the ability to quantify the assessment items.
\subsubsection{Virtual Technology Sets}\label{sec:2.5.1}
\indent The formulation and definitions of the \textit{virtual technology set} in the context of Step II within the PT and bTSc models are presented in  Eq. (\ref{Eq:41}) and Eq. (\ref{Eq:42}). Step I has similar expressions. In Eq. (\ref{Eq:41}) and Eq. (\ref{Eq:42}), the virtual scales of each $DMU_j$ (\textit{vInput}, \textit{vOuput}) and (\textit{avInput}, \textit{avOutput}) are defined in Eq. (\ref{Eq:12}) and Eq. (\ref{Eq:28}). The virtual scales in Step I for each $DMU_j$ represent certain aspects or characteristics within the model. These sets likely encapsulate specific parameters, constraints, or variables pertinent to the subsequent steps' formulation and computation. Virtual scales in Step II might represent refined or adjusted versions of the initial sets of Step I, potentially incorporating the outcomes or adjustments derived from the earlier steps.

\indent The transition from Step I to Step II usually involves refining or re-calibrating the parameters or sets based on the intermediate solutions or computations from the preceding step. This iterative process often helps converge toward more accurate or optimal results within the models.
\begin{equation}
\begin{aligned}
	\Phi_o^ {PT\star}&=\lbrace (\alpha^\star, \beta^\star) \mid (\alpha_j^{\star},\beta_j^{\star}), \forall j\in J\rbrace. 
	\end{aligned} 
	\label{Eq:41}
	\end{equation}
 \vspace{-1.5em}
\begin{equation} \begin{aligned} \label{Eq:42}  
 \Phi_o^ {bTSc\star}=
\lbrace (\alpha^{c\star},\beta^{c\star})\mid(\alpha_j^{c\star}, \beta_j^{c\star} ), \forall j\in J& \rbrace. 
	 \end{aligned}
	 \end{equation}	 
\subsubsection{Return-to-virtual-scale (RTvS)}\label{sec:2.5.2} Eq. (\ref{Eq:43}) and Eq. (\ref{Eq:44}) depict the computation steps for \textit{benchmark virtual scales} and \textit{affected benchmark virtual scale} in the PT and bTSc models during Step II  ( bvInput, bvOutput) =  ($\widehat\alpha_o^\star$, $\widehat\beta_o^\star$) and ( abvInput,  abvOutput) =  ($\widehat\alpha_o^{c\star}$, $\widehat\beta_o^{c\star}$), respectively. These equations likely involve transforming or deriving benchmark values based on specific conditions or constraints within the models.
\footnote{\label{fn:19}The virtual scale of each DMU is represented as a point on the 2D graphical intuition used to evaluate $DMU_o$. For example, Figures \ref{fig:3} and \ref{fig:4} show that the DMUs are located either beneath or on the b-Efficiency Equator when evaluating $DMU_K$. DEA models, however, cannot effectively express the PPS and its envelope. As a result, the points corresponding to efficient DMUs, as well as the projection point of $DMU_o$, may not lie on \textit{b-Efficiency Equator}.}

\indent The expressions for \textit{bvInput, bvOutput, abvInput, and abvOutput} in Step I are similar to those in  Eq. (\ref{Eq:43}) and Eq. (\ref{Eq:44}) with the replacement of the superscript "$\star$" with "$\#$." This step-wise process of deriving and refining benchmarks could indicate a progressive approach where the initial estimates from Step I are further modified or calibrated in Step II based on additional considerations or parameters introduced in the models.
\begin{equation} \label{Eq:43} 
	\widehat{\alpha}_o^\star= \sum_{\forall i\in I} \widehat{x}_{io} v_{io}^\star \textrm{  and  } \widehat{\beta}_o^\star= \sum_{\forall r\in R}\widehat{y}_{ro} u_{ro}^\star.	
	\end{equation}	
  \vspace{-1.5em}
\begin{equation}\label{Eq:44} \begin{aligned}	
	\widehat\alpha_o^{c\star} = [v_o^\star \widehat {x}_o + (1 - \gamma_o)\omega_o^{c\star} ]; \widehat\beta_o^{c\star} = (u_o^\star  \widehat{y}_o - \gamma_o\omega_o^{c\star} ).
\end{aligned}
\end{equation}
\indent  Eq. (\ref{Eq:45}) and Eq. (\ref{Eq:46}) compute the RTvS,  ($\Xi_o^ {PT\star}$ and  $\Xi_o^ {bTSc\star}$), is based on the virtual scales. The RTvS values might represent a measure or index that assesses the efficiency or productivity performance  $DMU_o$ within the context of these specific models. 
\footnote{\label{fn:20} The concepts of scale elasticity and returns to scale in the nonparametric methodology of DEA are crucial characteristics of a production function \cite{28}. The authors of DEA models assume the existence of a linear production frontier for $DMU_o$. In contrast, the bTSc and sTSc VGA models do not focus on the production function of the PPS. Instead, they measure the return-to-virtual-scale of $DMU_o$.}
\begin{equation}
	\Xi_o^ {PT\star}  = (\widehat\beta_o^\star/\beta_o^\star ) / (\widehat\alpha_o^\star/\alpha_o^\star)=1/E_o^ {PT\star}.
	\label{Eq:45}
	\end{equation}
  \vspace{-1em}
\begin{equation} 
 \label{Eq:46}
	\Xi_o^ {bTSc\star} = (\widehat\beta_o^{c\star} /\beta_o^{c\star}) / (\widehat\alpha_o^{c\star} /\alpha_o^{c\star} )
  =1/E_o^ {bTSc\star}.
\end{equation}

\indent For the case that $\omega_o^\star\geq \$0$, $\Xi_o^ {bTSc\star}$ and $E_o^ {bTSc\star}$ will be decreased and increased as the SIC scalar $\kappa_o^c$ rises. Conversely, when $\omega_o^\star\leq \$0$, $\Xi_o^ {bTSc\star}$ and $E_o^ {bTSc\star}$ will be decreased and increased, respectively, as the SIC scalar $\kappa_o^c$ diminishes. A preference is given to higher $\Xi_o^ {bTSc\star}$ and $E_o^ {bTSc\star}$. Depending on whether $\omega_o^\star\geq \$0$ or $\omega_o^\star\leq \$0$, $DMU_o$ has the flexibility to adjust $\kappa_o^c$ towards its lower or upper bound. This adjustment aims to achieve the final $\kappa_o^z$  that encompasses achievable benchmarks for inputs and outputs. The implications of these dynamics are further elucidated through the numerical examples and graphical illustrations presented in Section \ref{sec:5.1}, which depict how varying $\kappa_o^c$  influences the efficiency and the position of $DMU_o$ relative to the b-Efficiency Equator.
\footnote{\label{fn:21} Tables \ref{table:2} and \ref{table:3} listed the returns-to-virtual-scale of $DMU_o$. DEA models could not compute the return-to-scale, which is a concept that arises in the context of firms' production function in economics. Depending on the sign of the estimated $w_o^\star$, $DMU_o$ exhibits decreasing or increasing return-to-scale. }
\subsubsection{Interconnections Between Inputs and Outputs Indices} \label{sec:2.5.3}  Eq. (\ref{Eq:47}) and Eq. (\ref{Eq:48}) demonstrate adjustments to the affected virtual prices of an input $(v_{io}^\star x_{io}  + \gamma_{io}^Q \omega_o^{c\star} )$, and an output $(u_{ro}^\star y_{ro} - \gamma_{ro}^P \omega_o^{c\star} )$, respectively, considering the proportions to the vScalar, represented by  $\gamma_{io}^Q$  and  $\gamma_{ro}^P$. 
\begin{equation}\begin{aligned}
	\Delta_o^{bTSc\star} 
	= \sum_{\forall i\in I} (v_{io}^\star x_{io} + \gamma_{io}^Q \omega_o^{c\star} )  - \sum_{\forall r\in R}(u_{ro}^\star y_{ro} - \gamma_{ro}^P \omega_o^{c\star} ).
\end{aligned}	\label{Eq:47}
	\end{equation}	
 \vspace{-1.5em}
\begin{equation} \begin{aligned}
	\gamma_{io} ^Q  &= (1 - \gamma_o)q_{io}^\star/\sum_{\forall i\in I} q_{io}^\star, \forall i \in I ;\\
 \gamma_{ro}^P &= \gamma_o p_{ro} ^\star / \sum_{\forall r\in R} p_{ro}^\star, \forall r \in R.
	\end{aligned}\label{Eq:48}
	\end{equation}
\indent In Eq. (\ref{Eq:47}) and Eq. (\ref{Eq:48}), the proportions involving ($v_{io}^\star x_{io} +  \gamma_{io}^Q \omega_o^{c\star}$ ) and ($u_{ro}^\star y_{ro}-\gamma_{ro}^P \omega_o^{c\star} $) indeed highlight the interconnectedness between the input and output indices. These proportions indicate the adjustments made to the affected virtual prices concerning the quantities represented by $\gamma_{io}^Q$  and $\gamma_{ro}^P$, display how changes in one influence the other within the context of the model variables and constraints.

\subsubsection {2D Graphic Intuitions}\label{sec:2.5.4}
DMU-o can visualize the outcomes of the post-analysis based on the optimal solutions from the VGA models. For instance, we use the example dataset in Table \ref{table:1}. One can identify the third SIC scalar $\kappa_o^3$ by try and error. While the  $\kappa_o^3$ may not be within the range of  $\kappa_o^1$ and  $\kappa_o^2$. The particular $\kappa_o^3$ within the TS3 model produces $E_o^ {PT\star}=E_o^{TS3\star}$. Pure technical efficiency is affected by the SIC. 
The exact evaluations of PT and TS3 models are demonstrated in Figure \ref{fig:3}, in which the \textit{b-Efficiency Equator} should not be viewed as the efficiency frontier.

The \textit{anchor point}  
\footnote{\label{fn:22} Banker et al.\cite{29} used Figure \ref{fig:3} to illustrate the specific VRS PPS  with one input (x) and one output (y). The line segments passing through points A and E represent increasing and decreasing returns to scale, as indicated by the sign of $w_o^\star$. Each line intersects the vertical and horizontal axes at a distance of -$w_o^\star$, measured in the same units as x and y. For PPS with multiple inputs and outputs, $w_o^\star$ becomes dimensionless. See footnote \ref{fn:6}, the multiplier program \cite{8} of the Additive VRS Models normalizes input and output values to ensure compatibility with the dimensionless $w_o^\star$, a dual variable associated with the convexity condition. Additional details can be found in footnote \ref{fn:13}.

\indent Linear programming-based DEA models define decision variables and coefficients as dimensionless. However, this process disrupts the fundamental conditions of problem formulations, making it impossible to verify DEA models through duality.

\indent Each DEA model could not demonstrate the value of $w_o^\star$ on a graph. In each bTSc and sTSc VGA model, the anchor point is located on one of the axes of the 2D graphical intuition with intercept $\omega_o^\star (\$)$.
The actual measurement of the returns to scale in practice relies on the concept of \textit{scale elasticity} and {scale efficiency}, see Page 31 in Sickles and Zelenyuk\cite{7}. However, no one in the DEA literature has presented a solution for multiple inputs and outputs. }
of the bTSc model is represented by \textit{APc}, where "c" denotes the choices 1, 2, and 3. Points Kp, K3, AP3, and O are located at ($\alpha_o^{\star},\beta_o^{\star}$), ($\alpha_o^{3\star}, \beta_o^{3\star} $), ($[1-\gamma_o] \omega_ o^{3\star }, -\gamma_o  \omega_o^{3\star}$), and (0,0), respectively. Points O, K3, and AP3 are located at the corners of the triangle. The upper-right of the figure shows the two rectilinear distances between the pair points (Kp, Tp) and (K3, T3), the virtual gaps estimated by the PT and TS3 models. $DMU_o$ could comprehend the results between the PT and TS3 models to understand the effects of the SIC.

\indent Similarly, in Figure \ref{fig:4}, Points O, K1, and AP1 of the triangle are solutions of the TS1 model. The three points O, K2, and AP2 at the triangle are solutions of the TS2 model. Figure \ref{fig:4} showcases the effects of using different SIC scalars within a certain interval and helping select a final scalar $\kappa_o^z$ for productivity management. In this example, due to $w_o^\star >0$, we have $E_o^{TS2\star}$ $>$ $E_o^ {bTSz\star}$ $>$ $E_o^{TS1\star}$. The  bTSz model provides the most preferred $(q_o, p_o)$ ($q_o^ {bTSz\star}, p_o ^ {bTSz\star}$).

In Figure \ref{fig:4}, Points O, APc, and Kc at the triangle express solutions of the bTSc model. Eq. (\ref{Eq:49}), Eq. (\ref{Eq:50}), and Eq. (\ref{Eq:51}) are the expressions of the slopes for the vectors $\overrightarrow{\rm O, Kc}$, $\overrightarrow{\rm O, APc}$, and $\overrightarrow{\rm APc, Kc}$, respectively. The ratio of the scalar prices of inputs to outputs equals the slope $\bar{m}(\overrightarrow{\rm O, APc})$. Similarly, the ratio of the \textit{gray virtual input to output} equals the slope $\bar{m}(\overrightarrow{\rm APc, Kc})$.  
 \begin{equation} \begin{aligned}
	\bar{m}(\overrightarrow{\rm O,Kc})= \frac{{\beta_o^{c\star}-0}}{{\alpha_o^{c\star}-0}} =\frac{{\beta_o^{c\star}}}{{\alpha_o^{c\star}}}.
	\label{Eq:49}
	\end{aligned} \end{equation}
 \vspace{-1em}
\begin{equation}\begin{aligned}
	 \bar{m}(\overrightarrow{\rm O,APc})=\frac{-\gamma_o \omega_ o^{c\star }-0}{(1-\gamma_o) \omega_ o^{c\star }-0}=\frac{-\gamma_o \omega_ o^{c\star }}{(1-\gamma_o) \omega_ o^{c\star }}.
	\end{aligned}\label{Eq:50}
	\end{equation}	
 \vspace{-1em}
 \begin{equation}\begin{aligned}
	\bar{m}(\overrightarrow{\rm APc,Kc})=\frac{\beta_o^{c\star}-(-\gamma_o \omega_ o^{c\star })}{{\alpha_o^{c\star}}-[-(1-\gamma_o) \omega_ o^{c\star }]}
 =\frac{u_o^\star y_o}{ v_o^\star x_o}.
	\end{aligned}\label{Eq:51}
	\end{equation}	
 For instance, in Figure \ref{fig:3}, the triangle of the TS3 model has the following expressions. 
 \begin{equation} 
\begin{aligned}
	& \overrightarrow{\rm O,K3} =\overrightarrow{\rm O,AP3}+\overrightarrow{\rm AP3,K3} 
	 \iff \frac{\beta_o^{3\star}}{\alpha_o^{3\star}} =\frac{ u_o^\star y_o -\gamma_o  \omega_o^{3\star}}{ v_o^\star x_o +(1-\gamma_o)\omega_o^{3\star} }.
 \end{aligned} \label{Eq:52}
	 \end{equation}

\indent As shown in Eq. (\ref{Eq:53}), the estimated bTSc efficiency score could be visualized as the relationships of the three vectors.
\begin{equation} \begin{aligned}
 \overrightarrow{\rm O,Kc} =\overrightarrow{\rm O,APc}+\overrightarrow{\rm APc,Kc};\quad 
 E_o^ {bTSc\star}=\bar{m}(\overrightarrow{\rm O,Kc})= \frac{\beta_o^{c\star}}{\alpha_o^{c\star}}.
\label{Eq:53} \end{aligned}
\end{equation}
The effects of the SIC within the bTSc model could be visualized as triangular with points O, APc, and Kc in the 2D graphical intuition. The slopes $\bar{m}(\overrightarrow{\rm O,APc})$ and $\bar{m}(\overrightarrow{\rm APc, Kc})$ are not the scale and technical efficiencies. 

\section{Scenario-II: Evaluating the Super-Efficiency}\label{sec:3}
Figure \ref{fig:1} shows that Scenario II computes the super-efficiency. The PT and bTSc models are modified into sPT and sTEc models. When the efficient $DMU_o$ expands inputs and deduces outputs, its super-efficiency will be reduced to 1. In particular, $DMU_o$ is excluded from the reference set $J$ when assessing its super-efficiency. 
\subsection{super-Pure Technical Efficiency (sPT) Model}\label{sec:3.1}
\noindent \textbf{sPT model, TAP (dual) program:}
\footnote{\label{fn:23} Liu et al. \cite{30} reviewed slack-based measure of CRS efficiency and super-efficiency models. They attempted to tackle the existing problems in measuring the CRS super-efficiency. They offered two envelopment programs but without the multiplier programs. The SBM-based super-efficiency measurement models could not correctly evaluate serval $DMU_o$, see footnote \ref{fn:11} for the same reasons. Instead, using the sPT model could obtain the comprehensive assessments.}
\begin{equation}
	\delta_o^{sPT\star} = \min_{q_o, p_o, \pi_o} \sum_{\forall i\in I} q_{io} \tau_o  + \sum_{\forall r\in R}p_{ro} \tau_o;
	\label{Eq:54}
 \end{equation} \vspace{-1em}
\begin{equation}
	s.t. -\sum_{\forall j\in J-\{o\}} x_{ij} \pi_{jo} + q_{io} x_{io} \geq - x_{io}  , \forall i\in I;
 \label{Eq:55}
 \end{equation}  \vspace{-1.5em}
\begin{equation} 
 \qquad \sum_{\forall j \in J-\{o\}} y_{rj} \pi_{jo} + p_{ro} y_{ro} =  y_{ro}, \forall r \in R;	
\label{Eq:56}
 \end{equation}  \vspace{-2em}
\begin{equation} \label{Eq:57}
\qquad \pi_o, q_o  , p_o  \ge 0.	
 \end{equation}
\noindent \textbf{sPT model, TVG (primal) program:}	
\begin{equation} \label{Eq:58}
	\Delta_o^{sPT\star} = \max_{v_o, u_o} -\sum_{\forall i\in I} v_{io} x_{io}  + \sum_{\forall r\in R} u_{ro} y_{ro};
 \end{equation}  \vspace{-1em}
 \begin{equation} \begin{aligned} 
 \label{Eq:59}
	s.t. - \sum_{\forall i\in I} v_{io} x_{ij}  &+ \sum_{\forall r\in R}u_{ro}  y_{rj}  \le 0, \forall j\in J-\{o\}; 
 \end{aligned}
 \end{equation}  \vspace{-2em}
 \begin{equation}
	x_{io} v_{io} \ge\tau_o, \forall i\in I;	
 \label{Eq:60}
 \end{equation}  \vspace{-2 em}
\begin{equation}
	y_{ro} u_{ro}\ge\tau_o,\forall r\in R;
\label{Eq:61}
 \end{equation} \vspace{-2 em}
\begin{equation}
	v_o\geq 0\quad and \quad u_o \quad free.	
 \label{Eq:62}
 \end{equation}

Similar to Eq. (\ref{Eq:10}) and Eq. (\ref{Eq:11}), the optimal solutions of  Eq. (\ref{Eq:54}) and Eq. (\ref{Eq:58}) as depicted in Eq. (\ref{Eq:63}) and Eq. (\ref{Eq:64}). Usually, $\Delta_o^{sPT\#}$ in Step I exceeds \$1.  
\begin{equation}  \label{Eq:63} 
 \begin{aligned}
& \$0 <\delta_o^{sPT\#} = \sum_{\forall i\in I} q_{io} ^\# \times \$1  + \sum_{\forall r\in R}p_{ro} ^\#  \times \$1 =\Delta_o^{sPT\#} ; \\ 
&\$ 0 <\delta_o^{sPT\star } = \sum_{\forall i\in I} q_{io}^\star  \tau_o ^\star  + \sum_{\forall r\in R} p_{ro}^\star  \tau_o ^\star=\delta_{xo}^{sPT\star} + \delta_{yo}^{sPT\star} =\Delta_o^{sPT\star}< \$1.
  \end{aligned} 
 \end{equation} 
 \begin{equation}
\begin{aligned}
&\$0 < \Delta_j^{sPT\#}  =- v_o^\# x_j + u_o^\# y_j = -\alpha_j^\# + \beta_j^\# ;\\
&\$0 <\Delta_j^{sPT\star }  = -v_o^\star  x_j + u_o^\star  y_j =-\alpha_j^\star  + \beta_j^\star \leq \$1.
\end{aligned}
\label {Eq:64}
 \end{equation}

\indent  \textbf{DMU-o determines its unified goal price $\tau_o^\star$} to ensure $\delta_o^{sPT\star}$ and $=\Delta_o^{sPT\star}$ in Step II are within the range ($\$0, \$1$). Eq. (\ref{Eq:64}) supplied ($\alpha_o^\#, \beta_o^\#$) and ($\alpha_o^\star, \beta_o^\star$). 
 The relationship 
$\tau_o^\# :\tau_o^\star  =
1 :  \bar{t} = \Delta_o^{sPT\#}:\Delta_o^{sPT\star} $ 
could be derived from the \textit{relation equation} $\bar{t}   (  -\alpha_o^\#+ \beta_o^\# ) = 1(-\alpha_o^\star+\beta_o^\star$).  Use
 Eq. (\ref{Eq:65}) to obtain the dimensionless value of  $\bar{t}$. 
\begin{equation}
	\bar{t} = \$ 1/\beta_o^\# \textrm {  and  } \tau_o^\star   = \$ \bar{t}. 
	\label{Eq:65}
  \end{equation}
Dividing the relation equation by $\beta_o^\star$  substituting it $\bar{t}$, we have the new relation equation  $(1/\$\beta_o^{\star}) ( -\alpha_o^\#/\beta_o^\#+1) =( -\alpha_o^\star/\beta_o^\star+1)$. 

\indent Because of Eq. (\ref{Eq:58}) obtained $ \$ 0 \leq (-\alpha_o^\#+\beta_o^\# )$ and $\$0 \leq (-\alpha_o^\star+\beta_o^\star)$, we have $0 \leq ( -\alpha_o^\#/\beta_o^\#+1)$ and $0 \leq ( -\alpha_o^\star/\beta_o^\star+1)$. Using $\tau_o^\star=\$ \bar{t}$ in Step II, should have the solutions  $\beta_o^\star =\$1$ and $\alpha_o^\star<\$1$ to confirm $( -\alpha_o^\#/\beta_o^\#+1) =( -\alpha_o^\#/\beta_o^\#+1)$ and $ \$0 \leq \Delta_o^{sPT\star }$. The relative super-efficiency in Step I and Step II should be equal, as depicted in Eq. (\ref{Eq:66}). 
\begin{equation}
\begin{aligned} \label{Eq:66} 
	&1\leq E_o^{sPT} = \beta_o^\#  /\alpha_o^\#=\beta_o^\star  /\alpha_o^\star.	  
\end{aligned}
\end{equation}

 \indent \textbf{Numerical example.} The columns of B sPT-II and D sPT-II in Table \ref{table:3} summarized the evaluation solutions of $DMU_B$ and $DMU_D$ in Step II. The two efficient DMUs have different reference sets, ($DMU_A$) and ($DMU_B$, $DMU_G$). The super-efficiencies exceeded 1. The primal and dual programs of the sPT VGA model provide comprehensive evaluations. $DMU_B$ has the lowest value of $Y_1$ 567, and its super-efficiency $E_B^{sPT}$ equals 2.4126. $DMU_D$ has the highest value of $X_1$ 281, and its super-efficiency $E_D^{sPT}$ equals 1.3533.
 \footnote{\label{fn:24} DEA Additive Models could not measure the CRS efficiency; see footnotes \ref{fn:10} and \ref{fn:11} for the reasons. The existing Additive Models to measure the CRS super-efficiency do not verify the solutions by examining the dualities.}
 
 \indent Eq. (\ref{Eq:67}) likely represents this normalization process, where the Step I solutions are adjusted or transformed to achieve the solutions in Step II. 
\begin{equation} \begin{aligned}
  (\delta_o^{sPT\star },  \Delta_o^{sPT\star },  v_o^\star, u_o^\star)
  = \bar{t} \times (\delta_o^{sPT\#}, \Delta_o^{sPT\#}, v_o^\#, u_o^\#).
\label {Eq:67} \end{aligned}
\end{equation}
The symbol  $\mathcal{E}_o^{sPT}$  denotes the reference set of $DMU_o$ in the sPT evaluations. Each peer $DMU_j$ belonging to the reference set has sPT efficiency equal to 1 and its estimated intensity $\pi_{jo}^\star >0$. At the same time, the estimated intensity of the other inefficient DMUs, $\pi_{jo}^\star =0$. $DMU_o$ estimates the benchmark of each $X_i$ and $Y_r$,  $\widehat{x}_{io}^{sPT\star}$ and  $\widehat{y}_{ro}^{sPT\star}$ via  Eq. (\ref{Eq:68}). $DMU_o$ mimics the reference peers with their estimated intensities,  $\pi_{jo}^\star$, $\forall j\in \mathcal{E}_o^{sPT}$. $DMU_o$ is superior to its reference peers, whose criteria configurations are analog to $DMU_o$.
\begin{equation}
\begin{aligned}
	&\widehat{x}_{io}^{sPT\star} = \sum_{\forall j\in \mathcal{E}_o^{sPT}} x_{ij}  \pi_{jo}^\star  = x_{io} (1 + q_{io}^\star ), \forall i\in I; \\
	&\widehat{y}_{ro} ^{sPT\star} = \sum_{\forall j\in \mathcal{E}_o^{sPT}}y_{rj} \pi_{jo}^\star =y_{ro} (1 - p_{ro}^\star  )
 \forall r\in R.
\end{aligned}
\label {Eq:68}
\end{equation}
Assessing $DMU_o$ with the adjusted $X_i$ and $Y_r$, \(\widehat{x}_{io}^{sPT\star}\) and \(\widehat{y}_{ro} ^{sPT\star} \), as shown in Eq. (\ref{Eq:68}), will have the sPT efficiency   $\widehat{E}_o^{sPT }$ equals 1. 
\begin{equation}
    \sum_{\forall j\in \mathcal{E}_o^{sPT}}\pi_{jo}^\star   = \kappa_o^1 . 
 \label {Eq:69} \end{equation}
Let the total estimated intensities be the first SIC scalar $\kappa_o^1$, computed via  Eq. (\ref{Eq:69}). 

\subsection{Super-technical and Scalar choice Efficiency (sTSc) Model} \footnote{\label{fn:25} Chen et al.\cite{31} provided comprehensive reviews of VRS super-efficiency SBM-based models. They presented the two envelopment programs without the multiplier programs. Footnote \ref{fn:11} discussed the same reasons the two envelopment programs have incomplete evaluations. Instead, using $\kappa_o^{sTSc}=1$ to the sTSc model could obtain the comprehensive assessments.}
\label{sec:3.2}
  Adding the SIC Eq. (\ref{Eq:73}) to the sPT-TAP program will have the sTSc-TAP program. The SIC corresponds to the free-in-sign decision variable $w_o$. Choose the SIC scalar $\kappa_o^c$ to assess its super-efficiency $DMU_o$.\\
  
\noindent \textbf{sTSc model, TAP (dual) program:}

\begin{equation}\label {Eq:70} 
	\delta_o^{sTSc\star}= \min_{q_o, p_o} \sum_{\forall i\in I} q_{io} \tau_o + \sum_{\forall r\in R}p_{ro}   \tau_o;
\end{equation}  \vspace{- 1.5 em}
 \begin{equation} \label{Eq:71}
	s.t. -\sum_{\forall j\in J-\{o\}} x_{ij} \pi_{jo} + q_{io} x_{io}  \geq - x_{io}, \forall i \in I;	
 \end{equation}  \vspace{- 1.5 em}
\begin{equation}\label {Eq:72}
	 \sum_{\forall j\in J-\{o\}} y_{rj} \pi_{jo} + p_{ro} y_{ro} =  y_{ro}  , \forall r \in R;	
\end{equation}  \vspace{- 1.5 em}
 \begin{equation}\label {Eq:73}
	 \sum_{\forall j\in J-\{o\} } \pi_{jo}=\kappa_o^c;
  \end{equation}  \vspace{- 1.5 em}
  \begin{equation}\label {Eq:74}
	\pi_o,  q_o, p_o  \ge 0.	
\end{equation}
\noindent \textbf{sTSc model, TVG (primal) program:}	
\begin{equation}\begin{aligned}
	\Delta_o^{sTSc\star} =  \max_{x_o, v_o, w_o}  -\sum_{\forall i\in I} v_{io} x_{io}  &+ \sum_{\forall r\in R} u_{ro} y_{ro}  + \kappa_o^c  w_o;
 \label {Eq:75}
 \end{aligned}
 \end{equation}
  \vspace{- 1.5 em}
 \begin{equation} \begin{aligned}\label {Eq:76}
s.t. -\sum_{\forall i\in I} v_{io} x_{ij} + \sum_{\forall r\in R} u_{ro} y_{rj}  + 1  w_o  \le 0,
\forall j \in J-\{o\}; 
\end{aligned} \end{equation}  \vspace{- 2 em}
 \begin{equation}\label{Eq:77}
x_{io} v_{io} \ge\tau_o,\forall i \in I; 
\end{equation}  \vspace{-2em}
 \begin{equation}\label {Eq:78}
y_{ro} u_{ro} \ge \tau_o, \forall r \in R;	
 \end{equation}  \vspace{-2em}
 \begin{equation}\label {Eq:79}
v_o\geq 0, u_o,  w_o\quad free.
\end{equation}

\indent Similar to the bTSc model, the sTSc model determines the unified goal price in two steps.  \(\tau_o^\# = \$ 1\) and \(\tau_o^\star   = \$ \bar{t}\) will solve the TAP program comprehensively with the solutions \((q_o^\star , p_o^\star , \pi_o^\star )=(q_o^\#, p_o^\#, \pi_o^\# )\). Use Eq. (\ref{Eq:27}) to obtain \(\gamma_o\) and \((1- \gamma_o)\). Partition the \textit{vScalar} into two parts to reflect the effects of the SIC on inputs and outputs $DMU_o$.

\indent The maximized $\Delta_o^{sTSc}$ of $DMU_o$ in Step II of Eq. (\ref{Eq:75}) is expressed as the \textit{affected virtual Output (avOutput)} minus \textit{affected virtual Input (avInput)}, $-\alpha_o^{c\star}+ \beta_o^{c\star}$, as shown in Eq. (\ref{Eq:80}). The solutions in Step I have a similar expression, $-\alpha_o^{c\#}+ \beta_o^{c\#}$. Usually, in Step I $\Delta_o^{sTSc\#}$ exceeds \$1. 
\begin{equation}\begin{aligned} \label{Eq:80}
&\$0 \leq \Delta_o^{sTSc\star} = -v_o^\star x_o+u_o^\star y_o+\omega_o^{c\star }\\
&= -[v_o^\star x_o-(1-\gamma_o)\omega_o^{c\star }]+(u_o^\star   y_o+\gamma_o\omega_o^{c\star }) \\
& =-[gvInput-ivScalar]+[gvOutput+ovScalar]\\
	 &= -avInput+avOutput= -\alpha_o^{c\star} + \beta_o^{c\star}\leq \$ 1;\\
& \$0 \leq \Delta_o^{sTSc\#} =-\alpha_o^{c\#} + \beta_o^{c\#}.
 \end{aligned}  \end{equation}

 \indent Eq. (\ref{Eq:81}) expresses the solution of Eq. (\ref{Eq:76}) of each remaining $DMU_j$  $-\alpha_j^{c\star } + \beta_j^{c\star }$ in Step II. Similarly, the solution of Step I is  $-\alpha_j^{c\# } + \beta_j^{c\# }$. The estimated \textit{vScalar} of Step II  in  Eq. (\ref{Eq:76}), \(1w_o^*\), is decomposed into two components.
  \begin{equation}\begin{aligned} \label{Eq:81}
&\$0 < \Delta_j^{sTSc\star}=- v_o^\star x_j+u_o^\star y_j+ 1 w_o^{\star} \\
&= -[v_o^\star x_j-(1-\gamma_o) 1 w_o^{\star}]+(u_o^\star   y_j+\gamma_o 1 w_o^{\star}) = -\alpha_j^{c\star } + \beta_j^{c\star }, \forall j \in J,  j\neq {o} .
 \end{aligned} \end{equation}
Each $DMU_j$ is expressed by the pair of \emph{(avInput, avOutput)}, called as  \textit{virtual scales}. 

\indent  \textbf{DMU-o determines its unified goal price, $\tau_o^\star$} which is to ensure  $\delta_o^{sSTc\star}$ and  $\Delta_o^{sSTc\star}$ , in Step II, falls into the range of ($\$0, \$1$).
Eq. (\ref{Eq:80}) supplied ($\alpha_o^{c\#}, \beta_o^{c\#})$ and ($\alpha_o^{c\star}, \beta_o^{c\star}$). The relationship
$\tau_o^\# :\tau_o^\star  =
 1 :  \bar{t} = \Delta_o^{sTSc\#}:\Delta_o^{sTSc\star} $ could be derived from the relation equation $\bar{t}   (-\alpha_o^{c\#}+\beta_o^{c\#} ) = 1(-\alpha_o^{c\star}+\beta_o^{c\star}$). Use
 Eq. (\ref{Eq:82}) to obtain the dimensionless value of  $\bar{t}$. 
\begin{equation} \begin{aligned} \label {Eq:82}
	\bar{t }= \$ 1 / \beta_o^{c\#}= \$ 1 / (u_o^\#  y_o + \gamma_o \omega_o^{c\#})
 =\tau_o^\star.
\end{aligned} \end{equation}
Dividing the relative equation by $\beta_o^{c\star}$ and substituting $\bar{t}$, we have the new relation equation $1/\beta_o^{c\star}$  $ ( - \alpha_o^{c\#}/\beta_o^{c\#}+1) =( - \alpha_o^{c\star}/\beta_o^{c\star}+1)$.

\indent Because of Eq. (\ref{Eq:80}) obtained  $\$0\leq(-\alpha_o^{c\#}+\beta_o^{c\#} )$ and $ \$0\leq(-\alpha_o^{c\star}+\beta_o^{c\star})\leq\$1.$
Using $\tau_o^\star= \$ \bar{t}$ Step II, should have the solutions  $\beta_o^{c\star} =\$1,$ and  $\alpha_o^{c\star}\leq \$1$ , to confirm  $ (-\alpha_o^{c\#}/\beta_o^{c\#}+1) =( - \alpha_o^{c\star}/\beta_o^{c\star}+1)$. Therefore, we have $ \$0 \leq \Delta_o^{sTSc\star}\leq \$1$. Eq. (\ref{Eq:83}) shows the relative super-efficiency $E_o^{sTSc\star}$.
\begin{equation}
\begin{aligned}
 0 < E_o^{sTSc} &=\frac{u_o^\star y_o+\gamma_o \omega_o^{c\star}}{-v_{o}^\star x_o+(1-\gamma_o)\omega_o^{c\star} } =\beta_o^{c\star}/\alpha_o^{c\star}.	\end{aligned}
 \label {Eq:83}
  \end{equation}

\indent Normalizing the solutions of Step I by the dimensionless value $\bar{t}$ would obtain the solutions of Step II, as shown in  Eq. (\ref{Eq:84}). 
\begin{equation}\begin{aligned}\label {Eq:84}
	(\Delta_o^{sTSc\star }, \delta_o^{sTSc\star }, v_o^\star , u_o^\star ,  w_o^{\star})
	=  \quad\bar{t}\times(\Delta_o^{sTSc\#},\delta_o^{sTSc\#} ,v_o^\#, u_o^\# , w_o^{\#}).
\end{aligned}
\end{equation}

 The symbol  \(\mathcal{E}_o^{sTSc}\) denotes the set of reference peers in evaluating $DMU_o$ by the sTSc model, where each reference peer $DMU_j$ has \(\pi_{jo}^\star >0\). The other inefficient $DMU_j$ has \(\pi_{jo}^\star =0\). In  Eq. (\ref{Eq:81}), any best $DMU_j$ belongs to \(\mathcal{E}_o^{sTSc}\) has $\alpha_j^{c\#} = \beta_j^{c\#}$ and  $\alpha_j^{c\star } = \beta_j^{c\star }$, while the other remaining DMUs have $\alpha_j^{c\# } > \beta_j^{c\# }$ and $\alpha_j^{c\star } > \beta_j^{c\star }$. 

\indent Use  Eq. (\ref{Eq:85}) to compute the benchmarks of $X_i$ and $Y_r$, \(\widehat{x}_{io}^{sTSc\star}\) , and  \(\widehat{y}_{ro}^{sTSc\star}\). $DMU_o$ imitates the reference peers with their estimated intensities \((\pi_{jo}^\star )\).

\indent $DMU_o$ uses the benchmarks, $\widehat{x}_{io}^{sTSc\star}$ and $\widehat{y}_{ro} ^{sTSc\star}$, will decrease the sTSc super-efficiency $\widehat{E}_o^{sTSc }$ to 1 while the $\delta_o^{sTSc\star}$ drops to $\$$0.
\begin{equation}
\begin{aligned}\label {Eq:85} 
	\widehat{x}_{io}^{sTSc\star} &= \sum_{\forall j\in \mathcal{E}_o^{sTSc}}x_{ij}\pi_{jo}^\star = x_{io}(1 + q_{io}^\star ), \forall i\in I;\\
	\widehat{y}_{ro} ^{sTSc\star} &= \sum_{\forall j\in \mathcal{E}_o^{sTSc}}y_{rj}\pi_{jo}^\star = y_{ro}(1 - p_{ro} ^\star ), \forall r\in R .   \end{aligned} \end{equation}
\indent \textbf{Numerical example.} In Table \ref{table:3}, solutions of the sPT, sTS1, sTS2, and sTSz models to measure the relative super-efficiencies of $DMU_B$ and $DMU_D$ are summarized to observe how the decision variables are changed. $DMU_B$ and $DMU_D$ have the highest and lowest values on inputs and outputs and computed the super-efficiency scores. The four-phase procedure determines the range of SIC scalars ($\kappa_o^{sTS1}, \kappa_o^{sTS2})$. $DMU_B$ has the lowest value of $Y_1$ 567 and $DMU_D$ has the highest and lowest values of $X_2$ and $Y_2$, 281 and 97. The measured super-efficiencies are comprehensively measured due to the range of SIC scalars, which are predetermined in Phases 1 and 3. 
\footnote{\label{fn:26} Each of the sPT, sTS1, sTS2, and sTSz models determines $\tau_o^\star$ in two steps. Those VGA models have comprehensive evaluations, but DEA models cannot.  }

\subsection{The SIC Effects in Measuring the Super-Efficiency}\label{sec:3.3}
In Figure \ref{fig:5}, the coordinates of points O, AP1, AP2, B1, and B2 of the sTS1 and sTS2 models in evaluating $DMU_B$ can be visualized. Eq. (\ref{Eq:86}), Eq. (\ref{Eq:87}), and Eq. (\ref{Eq:88}) compute the slopes of the vectors $\overrightarrow{\rm O, Bc}$, $\overrightarrow{\rm O, APc}$, and $\overrightarrow{\rm APc, Bc}$, respectively, where the choice "c" could be "1" and "2."  
\begin{equation} \begin{aligned}
	\bar{m}(\overrightarrow{\rm O,Bc})= \frac{{\beta_o^{c\star}-0}}{{\alpha_o^{c\star}-0}} =\frac{{\beta_o^{c\star}}}{{\alpha_o^{c\star}}}.
	\label{Eq:86}
	\end{aligned} \end{equation}	
 \vspace{-1.5em}
\begin{equation}\begin{aligned}
	 \bar{m}(\overrightarrow{\rm O,APc})=\frac{\gamma_o \omega_ o^{c\star }-0}{-(1-\gamma_o) \omega_ o^{c\star }-0}=\frac{\gamma_o \omega_ o^{c\star }}{-(1-\gamma_o) \omega_ o^{c\star }}.
	\end{aligned}\label{Eq:87}
	\end{equation}	
 \vspace{-1em}
 \begin{equation}\begin{aligned}
	\bar{m}(\overrightarrow{\rm APc,Bc})=\frac{\beta_o^{c\star}-\gamma_o \omega_ o^{c\star }}{{\alpha_o^{c\star}}-[-(1-\gamma_o) \omega_ o^{c\star }]}=\frac{u_o^\star y_o}{v_o^\star x_o}.
	\end{aligned}\label{Eq:88}
	\end{equation}	
 The same analysis should apply to $DMU_D$ the situation as depicted in Figure \ref{fig:5} 5.

\indent The three vectors have the relationships as shown in Eq. (\ref{Eq:89}). The estimated sTSc super-efficiency score is affected by the SIC. 
\begin{equation} \begin{aligned}
 \overrightarrow{\rm O,Bc} =\overrightarrow{\rm O,APc}+\overrightarrow{\rm APc,Bc};\quad
 E_o^{sTSc\star}=\bar{m}(\overrightarrow{\rm O,Bc})= \frac{\beta_o^{c\star}}{\alpha_o^{c\star}}.
\label{Eq:89} \end{aligned}
\end{equation}
The effects of the SIC within the sTSc model could be visualized as triangular with points O, APc, and Bc in the 2D graphical intuition. The slopes $\bar{m}(\overrightarrow{\rm O,APc})$ and $\bar{m}(\overrightarrow{\rm APc, Bc})$ are not the scale and technical efficiencies. 

\section{Four-phase b-VGA-EA method}\label{sec:4}
The following presentation is to calibrate $(q_o, p_o)$ to measure the efficiency of $DMU_o$. It can be applied to measure the super-efficiency. Nonetheless, the PT model estimate $(q_o, p_o)$ may not be achievable $DMU_o$. The bTSc model with the SIC equates to the choice scalar $\kappa_o^c$ that identifies the reference peers $\mathcal{E}_o^ {bTSc}$. The benchmarks of inputs and outputs are the function of intensities of the reference peers belonging to $\mathcal{E}_o^ {bTSc}$, as shown in Eq. (\ref{Eq:33}). Section \ref{sec:2.5.4} illustrates the effects of SIC on the solutions.    

\indent The four-phase process illustrated in Figure \ref{fig:1} significantly improves the EA of DMUs. This framework is structured around conditions that epitomize a performance improvement problem. It incorporates practical considerations into the PT and bTSc VGA models for systematic solutions. 
\subsection{Phase-1: PT model Identifies the First SIC Scalar} \label{sec:4.1} 
An initial assessment distinguishes between efficient and inefficient DMUs in the performance evaluation. This distinction is crucial for the subsequent analysis and is primarily based on the existence of a "virtual gap." An inefficient $DMU_o$ is characterized by a positive virtual gap, indicating a discrepancy between its current performance and the b-Efficiency Equator. The DMU must undergo a comprehensive four-phase improvement process.
  
\indent The commencement of this process involves the identification of the first SIC scalar $\kappa_o^1$. This scalar is conceptualized as the sum of the estimated intensities associated with efficient counterparts $\kappa_o^1 = \sum_{\forall j\in \mathcal{E}_o^ {PT}}\pi_{jo}^\star$. Essentially, it serves as a quantitative measure of the extent to which $DMU_o$ needs to adjust its operations to align with the practices of its reference peers.

\indent Conversely, a DMU identified as efficient exhibits a zero virtual gap, signifying that its performance already aligns with the best practices observed across the dataset at the optimal level. Such DMUs are deemed ineligible for further analysis in Scenario II. For these units, Scenario II presented in Section \ref{sec:3} should be employed for continuous performance management and improvement, ensuring they maintain their efficiency status over time.

\subsection {Phase-2: TS1 Model Uses the First SIC Scalar}\label{sec:4.2}
The comparison between the optimal solutions of the PT and TS1 models in Steps I and II  reveals some critical insights. Both models share identical TAP solutions in Step I due to the equality $\tau_o^ {PT\#} = \tau_o^{TS1\#}$ and $\sum_{\forall j\in \mathcal{E}_o^ {PT}} \pi_{jo}^\#  $= $\sum_{\forall j\in \mathcal{E}_o^{TS1} } \pi_{jo}^\#$ = $\kappa_o^1$. However, their TVG programs yield two sets of optimal solutions due to the additional \textit{vScalar}. This process leads to a specific relationship between their intensities, resulting in disparate benchmarks for performance indices between the two models.

\indent Eq. (\ref{Eq:90}) shows Step I of the PT and TS1 models have the relationships.
\begin{equation} \label{Eq:90}
\begin{aligned}
	&\Delta_o^{TS1\#} = \Delta_o^{PT\#}= \delta_o^{TS1\#} = \Delta_o^{PT\#} ; (q_o^ {PT\#}, p_o ^ {PT\#})=(q_o^{TS1\#}, p_o^{TS1\#});\\
	& (v_o^ {PT\#}, u_o^ {PT\#} )\neq(v_o^{TS1\#}, u_o^{TS1\#}); \sum_{\forall j\in \mathcal{E}_o^ {PT}} \pi_{jo}^\# = \sum_{\forall j\in \mathcal{E}_o^{TS1}} \pi_{jo}^\# = \kappa_o^1.   \end{aligned} \end{equation}	
    
 Note the best peers’ intensities, $\pi_{jo}^ {PT\#}$, and $\pi_{jo}^{TS1\#}$, may be different; therefore, the benchmarks of performance indices of the two models are distant. They have identical solutions on the dual variables. But have different primal solutions. 

\indent Step II further differentiates the two models. Although their TAP programs provide the same total values, their goal prices and other variables differ. This disparity leads to different efficiency measurements: $E_o^ {PT\star}$  represents pure technical efficiency, and $E_o^{TS1\star}$ comprises the effects on the \textit{gray virtual input and output} and scalar prices. Notably, the SIC scalar is a crucial factor in evaluating $DMU_o$. In Step II of PT and TS1 models,  the TAP programs equal to
($\sum_{\forall i\in I}$ $q_{io}^\star$ 
+ $\sum_{\forall r\in R}$ $p_{ro}^\star$ ) multiplies the distinct virtual goal prices $\tau_o^ {PT\star}$  $\tau_o^{TS1\star}$. The PT and TS1 models have the relationships shown in  Eq. (\ref{Eq:91}).
\begin{equation}\label{Eq:91}
\begin{aligned}
	&\tau_o^{TS1\star} \neq \tau_o ^ {PT\star}, \Delta_o^{TS1\star} \neq \Delta_o^{PT\star};\quad(v_o^ {PT\star}, u_o^ {PT\star}) \neq (v_o^{TS1\star}, u_o ^{TS1});\\ 
	&(q_o^ {PT\star}, p_o ^ {PT\star}) = (q_o^{TS1\star}, p_o ^{TS1\star});\quad 0 < E_o^ {PT\star} \neq E_o^{TS1\star} \le1.	
	\end{aligned}
	\end{equation}

\subsection {Phase-3: Identifying the Second SIC Scalar} \label{sec:4.3} The standard sensitivity analysis of the LP TS1 model involves perturbing the $\kappa_o^1$  scalar, which results in allowable decreases and increases. Specifically, $\kappa_o^2$  equals $\kappa_o^1$  minus the allowable decrease or $\kappa_o^2$  equals $\kappa_o^1$  plus the allowable increase, as expressed in  Eq. (\ref{Eq:92}).
\begin{equation} \begin{aligned}
	\kappa_o^2  &= \kappa_o^1  - \textrm{(allowable decreasing of } \kappa_o^1), \\
 \textrm{or  }
	 \kappa_o^2   &= \kappa_o^1   + \textrm{(allowable increasing of }\kappa_o^1 ).
	\end{aligned} 
	\label {Eq:92}
	\end{equation}	
 
The final SIC scalar $\kappa_o^c$ in Phase-4 is linked to $\kappa_o^1$ in Phase-1. These bounds ($\kappa_o^1$  and $\kappa_o^2$) represent the possible range of SIC scalars ($\kappa_o^c $) within the bTSc model. If $w_o^{TS1\star}$  and $w_o^{TS2\star}$ are more significant than 0, $\kappa_o^1  < \kappa_o^c  < \kappa_o^2 $. Otherwise, $\kappa_o^1  > \kappa_o^c  > \kappa_o^2 $. 

\indent The best peers in the TS1 and  TS2 models remain consistent, as indicated by $\mathcal{E}_o^{TS1}$=$\mathcal{E}_o^{TS2}$=$\mathcal{E}_o^ {bTSc}$. However, the estimated intensities of these peers differ $\pi_o^{TS1\star}$ $\neq\pi_o^{TS2\star}$. The solutions of these models exhibit relationships detailed in Eq. (\ref{Eq:93}), displaying discrepancies in the estimated variables, such as:
\begin{equation}
\begin{aligned}
 \tau_o ^{TS1\star} \neq \tau_o ^{TS2\star};
\quad(v_o^{TS1\star}, u_o^{TS1\star}) \neq (v_o^{TS2\star}, u_o^{TS2\star});\\  w_o^{TS1\star}\neq w_o^{TS2\star};\quad (q_o^{TS1\star}, p_o^{TS1\star}) \neq(q_o^{TS2\star}, p_o^{TS2\star}).
 \label{Eq:93}  
 \end{aligned} 
 \end{equation} 

\subsection{Phase-4: Select the Final Scalar Within the Region}\label{sec:4.4} 
Phases -1, -2, and -3 integration confirms $DMU_o$'s learning process with the reference peers, whose total intensities are bounded within the $\kappa_o^1$  and $\kappa_o^2$. Phase-4, a critical decision-making step, involves selecting the final scalar $\kappa_o^z$  for the   bTSz model. $DMU_o$ undertakes several trials in choosing scalars in Phase-4 to determine a preferred final scalar $\kappa_o^z$. Each choice scalar corresponds to a unique set of estimated intensities and benchmarks of inputs and outputs  $\widehat{x}_{io}^ {bTSz}$  and  $\widehat{y}_{ro}^ {bTSz}$ derived from efficient peers.

\indent  By emulating the reference peers based on the estimated intensities, $DMU_o$ gains insights into feasible and desirable productivity benchmarks. This interactive approach ensures a more thorough evaluation $DMU_o$. Exploring different SIC scalars and model variations aids in comprehending the diverse impacts these factors have on assessments. Ultimately, this process helps make better-informed management decisions regarding productivity and efficiency.

DMU-o may redefine the datasets by altering the performance indices and DMUs and repeating the four-phase procedure. 

\section{Numerical examples and the 2D geometric intuitions}\label{sec:5} The numerical example depicted in Table \ref{table:1}. Table \ref{table:2} summarizes the efficiency measurements of $DMU_K$, $DMU_o$, within Scenario-I by the PT, TS1, TS2, and TS3 models. The PT model evaluated $DMU_B$ and $DMU_D$ and confirmed they are efficient units. Section \ref{sec:5.2} offers Scenario II to measure the super-efficiencies of efficient DMUs.
\subsection{Measures the Inefficiencies}\label{sec:5.1} 
\subsubsection{Relationships Between PT and TS1 Models in Step I}\label{sec:5.1.1}
Eq. (\ref{Eq:18}) and Eq. (\ref{Eq:35}) are used to calculate the total adjustment prices ($\Delta_o^{PT\#}$ = $\delta_o^{TS1\#}$) and total virtual gaps ($\Delta_o^{PT\#}$= $\Delta_o^{TS1\#}$), using specific data provided in row R2 of Table \ref{table:2}. Each amount is \$2.3010. Additionally, the TS1-I model incorporates the SIC impact through the $\textit{vScalar}$, $\omega_o^{1\#}$(=\$2.479), in calculating the total virtual gap. This result demonstrates how the virtual gap, considering the SIC and vScalar in the TS1-I model, maintains the same value as in the PT-I model, indicating a consistent evaluation despite including additional factors.
\subsubsection{Compare the PT and TS3 Models}\label{sec:5.1.2}
By try and error, we identified $\kappa_o^3$ equals 0.718 to have $E_o^ {PT\star}$=$E_o^{TS3\star}$(=0.589). The columns of PT-II and TS3-II in Table \ref{table:2} list the estimated solutions, in which rows R8 and R9 are the two virtual technology sets, ($\alpha_j^\star,\beta_j^\star$) and ($\alpha_j^{3\star},\beta_j^{3\star}$), where $DMU_j$= (A, B, D, G, H, K), as per  Eq. (\ref{Eq:12}),  Eq. (\ref{Eq:28}), and Eq. (\ref{Eq:29}). 

  \begin{figure}[ht]
\centering
\includegraphics[width=0.64\linewidth, height=10.5cm]{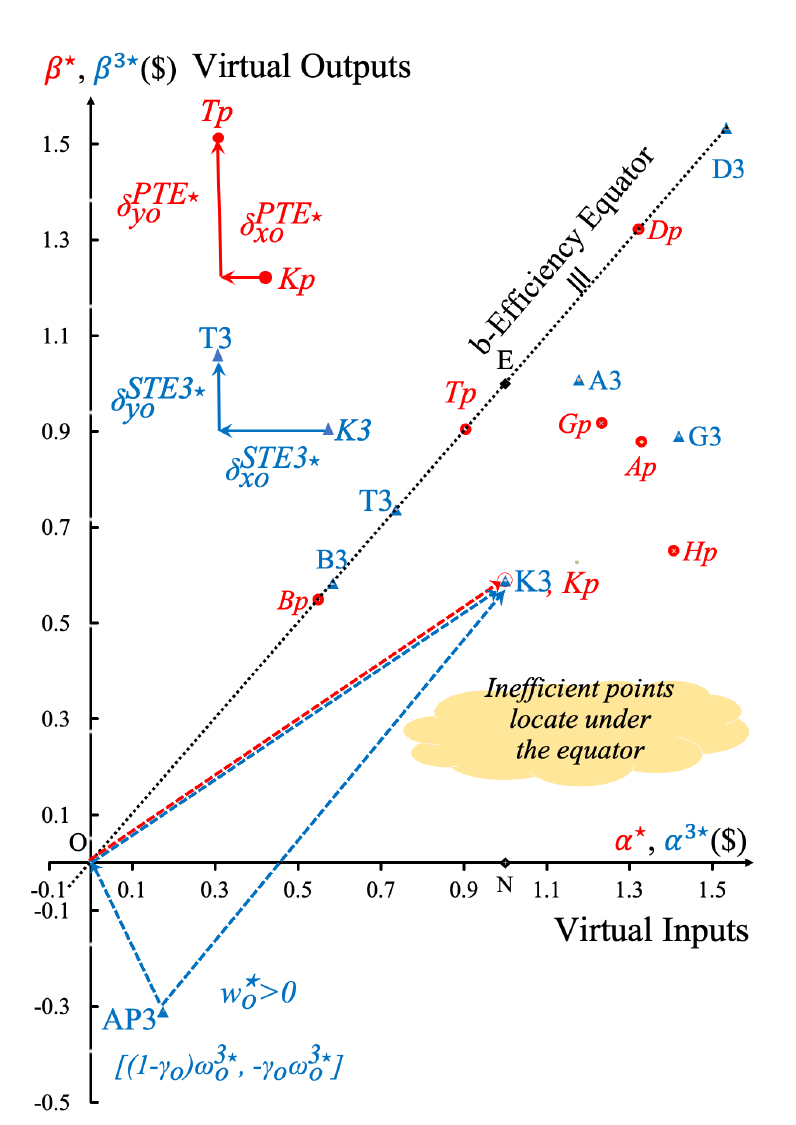} 
\caption{Comparing the solutions of PT and TS3 models.}
\label{fig:3}
\end{figure}
\indent Via  Eq. (\ref{Eq:43}) and Eq. (\ref{Eq:44}), the locations of Tp, ($\widehat\alpha_o^\star,\widehat\beta_o^\star$) and point T3, ($\widehat\alpha_o^{3\star},\widehat\beta_o^{3\star}$) on the b-Efficiency Equator are obtained. The estimated rectilinear distances from points Kp to Tp and from point K3 to T3 are $\Delta_o^{PT\star}$, and $\delta_o^{TS3\star}$, see the expressions in the upper-left corner of Figure \ref{fig:3}. Section \ref{sec:2.5.4} illustrated the SIC affected the estimations as shown in Eq. (\ref{Eq:52}). 

 \indent In Figure \ref{fig:3}, each point signifies a $DMU_j$ within the virtual technology sets of the PT-II ($\Phi_o^ {PT\star}$) and TS3-II models ($\Phi_o^{TS3\star}$), respectively. 
 The points (Ap, Bp, Dp, Gp, Hp, Kp) and (A3, B3, D3, G3, H3, K3) are located at ($\alpha_j^\star,\beta_j^\star$) and ($\alpha_j^{3\star},\beta_ j^{3\star}$). The blue dashed lines graphically represent the vectors $\overrightarrow{\rm O, K3}$ $\overrightarrow{\rm AP3, K3}$. The red dashed line indicates the vector $\overrightarrow{\rm O, Kp}$ of pure technical efficiency,  $E_o^ {PT\star}$= $\beta_o^{3\star}$/$\alpha_o^{3\star}$.  Inefficient DMUs are situated under the b-Efficiency Equator.

\subsubsection {Compare the Solutions between TS1 and TS2 models} \label{sec:5.1.3} 
\begin{figure}[ht]
\centering
\includegraphics[width=0.64\linewidth, height=9cm]{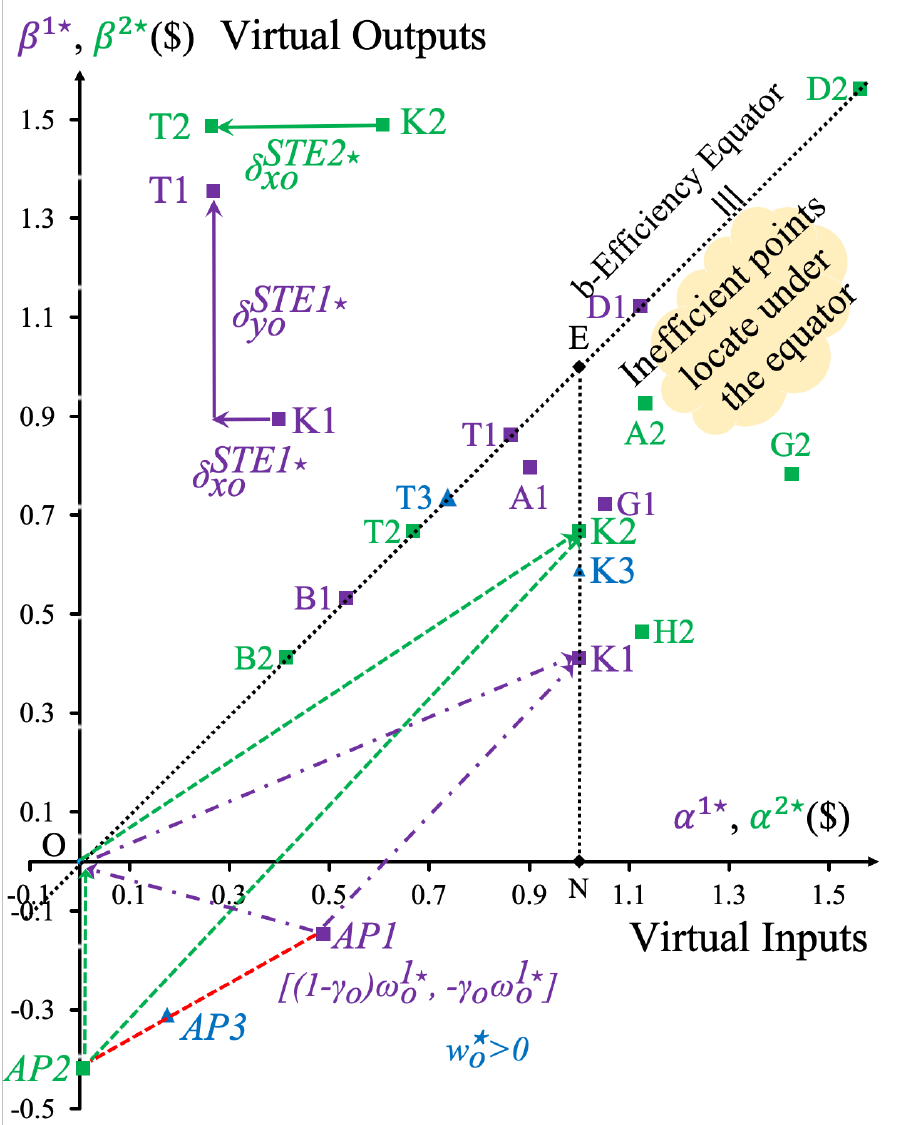}
\caption{Comparing the solutions of TS1 and TS2 models. \hfill}
\label{fig:4}
\end{figure}

Figure \ref{fig:4}, similar to Figure \ref{fig:3}, depicts the locations of DMUs in the TS1-II and TS2-II models, respectively, (A1, B1, D1, G1, H1, K1)  and (A2, B2, D2, G2, H2, K2). Their coordinates, ($\alpha_j^{1\star},\beta_j^{1\star}$) and ($\alpha_j^{2\star},\beta_j^{2\star}$), are listed in rows R8 and R9.  DMUs B and D emerge as efficient peers with their intensities $\pi_{oB}^\star$  and $\pi_{oD}^\star$.  R4 detailed the benchmarks of inputs and outputs. R5 lists the virtual prices. R6 contains the benchmark virtual scales. R7 lists $\Xi_o^ {bTSc\star}$ and $E_o ^ {bTSc\star}$. Since $w_o^\star>0$, decreasing $\kappa_o^1$  (=1.5153)  to  $\kappa_o^2$  (=0.515), $DMU_o$ exhibits decreasing RTvS: $\Xi_o^{TS1\star} (=2.435) > \Xi_o^{TS2\star} (=1.497)$.

\indent At the same time, the relative efficiency is increasing, $E_o^{TS1\star}$ (=0.411) less than $E_o^{TS2\star}$ (=0.688). However, $DMU_o$ may use  Eq. (\ref{Eq:48}) to analyze the interlinkage relationships between input and output indices for managing the performance indices.

\indent The bTSc program introduces an \textit{anchor point}  (APc) located at the point $[(1 -\gamma_o)\kappa_o^c  w_o^\star$, $- \gamma_o \kappa_o^c w_o^\star]$ denotes as [AP(x), AP(y)]. Anchor points are located in the second and fourth quadrants when $w_o^\star> 0$ and $w_o^\star< 0$. The rectilinear distance between point APc and the origin O is the vScalar, $\kappa_o^c w_o^\star$. In Figure \ref{fig:4}, the anchor points of the TS1 and TS2 models are located at  AP1 and  AP2. Because $w_o^\star>0$, the two anchor points are located in the fourth quadrant. The TS3 model interpolates TS1 and TS2 models, providing a framework for determining the optimal scalars $\kappa_o^3$ within the interval of SIC scalars. Using the scalar $\kappa_o^3$ between $\kappa_o^1$  and $\kappa_o^2$, the TS3 model has the anchor points, and AP3 is located between AP1 and AP2.  

 \indent  $DMU_j$ lies in the fourth quadrant if, according to  Eq. (\ref{Eq:2}), one could have $\Delta_j^ {bTSc\star}> \$0$. ($\alpha_j^{c\star}>\$0, \beta_j^{c\star}<\$0$), and their ratio is less than zero. Conversely, $DMU_j$ is in the third quadrant if ($\alpha_j^{c\star}<\$0, \beta_j^{c\star}<\$0$), and their ratio is more significant than zero. 

\indent Section \ref{sec:2.5.4} illustrated the SIC affected the estimations as shown in Eq. (\ref{Eq:53}). Therefore, eliminating an inefficient $DMU_j$ from (\textit{X, Y}) does not impact the evaluation. We can visually interpret the effects under the SIC scalars $\kappa_o^1$  and $\kappa_o^2$ using geometric vectors. In Figure \ref{fig:4}, the purple color dashed lines express the three vectors $\overrightarrow{\rm O,K1},  \overrightarrow{\rm AP1,O}$, and $\overrightarrow{\rm AP1,K1}$ of the TS1 model can be found.

 \indent  Similarly, the green color dashed lines express the three vectors $\overrightarrow{\rm O, K2},  \overrightarrow{\rm AP2, O}$, and $\overrightarrow{\rm AP2, K2}$ of the TS2 model. The three vectors $\overrightarrow{\rm O,K3}$ $\overrightarrow{\rm AP3,O}$  $\overrightarrow{\rm AP3, K3}$ of the TS3 model are not depicted. Point K3 is located between points K1 and K2 on the normalization vertical line (N, E).  Using the final SIC scalar $\kappa_o^z$ in the  bTSz model, points APz, Kz, and Tz are located between points (AP1 and AP2), (K1 and K2), and (T1 and T2), respectively. The three vectors, $\overrightarrow{\rm O,Kz}$, $\overrightarrow{\rm APz,O}$, $\overrightarrow{\rm APz,Kz}$, shall be visualized in Figure \ref{fig:4}, similar to Eq. (\ref{Eq:53}). The upper left corner of the figure depicts the rectilinear distances of $DMU_o$ the project on the b-Efficiency Equator. 
\begin{table}[b]
\centering
\setlength{\tabcolsep}{2pt}
\renewcommand{\arraystretch}{0.91}
\caption{Evaluate Inefficiencies for $DMU_K$, $DMU_B$, and $DMU_D$ by Scenario-I. } \label{table:2}
\footnotesize
\begin{tabularx}{\linewidth}{lccrrrrrrrrr} 
  & &  $DMU_o$= & $K$ & $K$ & $K$ & $K$ & $K$ & $K$ & $K$ & $B$ & $D$\\ 
 \hline
  Row& Sol. & Unit &PT-I&PT-II &TS1-I&TS1-II&TS2-I&TS2-II &TS3-II&PT-II&PT-II\\
  & & & \multicolumn{2}{c} {Phase-1} & \multicolumn{2}{c} {Phase-2}&\multicolumn{2}{c}{Phase-3}&{Ph.-4}&{Ph.-1}&{Ph.-1}\\
  \hline
R1 &$\kappa_o^c$   & - &1.5153 &1.5153 &1.5153 &1.5153 &0.5150 &0.5150& 0.718&1 &1\\
   &$\tau_o$ &  $\$$ &1 &0.179 &1.000 &0.256 &1.000 &0.500 &0.413&0.500&	0.266\\
\hline
R2 &$\Delta_o,\delta_o$ & \$ &2.3010 &0.4113 &2.3010 &0.5893 &0.6643 &0.3321&0.411 &0 &0 \\
   &$v_{1o}$  &  \$/ton &2.8713 &0.5133 &0.6250 &0.1601 &0.6250 &0.3125 &0.258&0.500	&0.387 \\
   &$v_{2o} $  &  \$/hr &0.0069 &0.0012 &0.0069 &0.0018 &0.0069 &0.0034&0.003 & 0.017	&0.001\\
   &$u_{1o}$ &  \$/$m^3$ &0.0022 &0.0004 &0.0011 &0.0003 &0.0011 &0.0006&0.0005&0.001&	0.0003 \\
   &$u_{2o}$  &  \$ / $\%$ &0.0204 &0.0036 &0.0204 &0.0052 &0.0204 &0.0102 &0.008&0.006&	0.003 \\
   &$w_o$ &  \$  &0 &0 &1.6362 &0.4190 &1.6362 &0.8181&0.675 &0&	0 \\
\hline
R3 &$q_{1o} $  & - &0 &0 &0 &0 &0.4554 &0.4554&0.363&0 &0 \\
   &$q_{2o}$   &-  &0.5334 &0.5334 &0.5334 &0.5334 &0.2089 &0.2089 &0.275&0 &0 \\
   &$p_{1o}$   &-  &0 &0 &0 &0 &0 &0&0&0 &0 \\
   &$p_{2o} $  &-&1.7677 &1.7677 &1.7677 &1.7677 &0 &0 &0.359&0 &0 \\
   & $\pi_{oB}$ &- &1.421 &1.421 &1.421 &1.421 &0.119 &0.119&0.383&1&0 \\
   & $\pi_{oD}$  &- &0.094 &0.094 &0.094 &0.094 &0.396 &0.396 &0.335&0 &1 \\
   \hline
R4 &$\widehat{x}_{1o}$ &ton &1.6 &1.6 &1.6 &1.6 &0.8713 &0.8713&1.019 &1.0&	1.9\\
   &$\widehat{x}_{2o}$  &hr &67.66 &67.66 &67.66 &67.66 &114.72 &114.72 &105.165&29	&281\\
   & $\widehat{y}_{1o}$  & $m^3$ &1036 &1036 &1036 &1036 &1036 &1036 &1036&567&	2446 \\
   & $\widehat{y}_{2o}$ &\%  &135.6 &135.6 &135.6 &135.6 &49.0 &49.0&66.58&89&	97\\
   & $\gamma_o$ & -  &0.232 &0.232 &0.232 &0.232 &1.000 &1.000&0.640&0&0 \\
   &  $\omega_o^c$ &  \$ &0 &0 &2.479 &0.635 &0.843 &0.421&0.485 &0 &0 \\
   \hline
R5 &$x_{1o}v_{1o}$ &  \$ &4.594 &0.821 &1.000 &0.256 &1.578 &0.789&0.689&0.500&	0.734 \\
   &$x_{2o}v_{2o}$ &  \$ &1.000 &0.179 &3.479 &0.891 &1.265 &0.632 &0.622&0.500&	0.266\\
   &$y_{1o}u_{1o}$ &  \$ &2.293 &0.410 &1.178 &0.302 &1.178 &0.589 &0.486&0.500	&0.734 \\
   &$y_{1o}u_{1o}$ &  \$ &1.000 &0.179 &3.479 &0.891 &1.843 &0.921 &0.898&0.500	&0.266 \\
\hline
R6 &$\alpha_o,\alpha_{o}^c$&  \$ &5.594 &1.000 &3.905 &1.000 &2.000 &1.000&1.000 &1 &1 \\
&$\beta_o,\beta_{o}^c$ &  \$ &3.293 &0.589 &1.604 &0.411 &1.336 &0.668&	0.589 &1 &1 \\
&$\widehat\alpha_o,\widehat\alpha_{o}^c$ &  \$ &5.061 &0.905 &3.371 &0.863 &1.336 &0.668&0.737 &1&1 \\
   &$\widehat\beta_o,\widehat\beta_{o}^c$ &  \$&5.061 &0.905 &3.371 &0.863 &1.336 &0.668 &0.737&1 &1\\
\hline
R7 & $\Xi_o$ &- &1.699 &1.699 &2.435 &2.435 &1.497 &1.497 &1.699&1&1 \\
   & $E_o$ & - &0.589 &0.589 &0.411 &0.411 &0.668 &0.668&0.589 &1 &1 \\
   \hline 
     R8&$\alpha_A,\alpha_{A}^c$&  \$ &7.432 &1.000 &3.522 &0.902 &2.265 &1.133&0.902&3.219	&1.002 \\ 
   &$\alpha_B,\alpha_{B}^c$&  \$ &3.071 &1.328 &2.082 &0.533 &0.825 &0.413 &0.533&1&	0.414 \\   
      &$\alpha_D,\alpha_{D}^c$&  \$ &7.393 &0.549 &4.382&1.122 &3.125 &1.563 &1.122&5.795&	1 \\ 
   \hline 
  R9&$\beta_A,\beta_{A}^c$ &  \$ &4.917 &0.879 &3.110 &0.796 &1.853 &0.926 &0.796&1.715&	0.664 \\
   &$\beta_B,\beta_{B}^c$ &  \$ &3.071 &0.549 &2.082 &0.533 &0.825 &0.413 &0.533&1&	0.414 \\
   &$\beta_D,\beta_{D}^c$ &  \$ &7.393 &1.322 &4.382 &1.122 &3.125 &1.563 &1.122&2.702&	1 \\
   \hline
  R10&AP(x)&\$& -&	-	&1.905	&0.488&	0&	0&	0.175&	-	&-\\
&AP(y)&\$&-&	-&	-0.575	&-0.147	&-0.843&	-0.421	&-0.310	&-	&-\\ \hline
 \end{tabularx}
 \end{table}

\subsection{Measures the Super-efficiencies}\label{sec:5.2} Scenario-I confirmed $DMU_B$ and $DMU_D$ are efficient.  Scenario II employs the sPT and sTSc models to measure their super-efficiencies in the four-phase procedure. Table \ref{table:3} summarizes the solutions of the sPT, sTS1, sTS2, and sTSz models. For simplicity, in rows R8 and R9, the values of $\alpha_{oK}^\star, \beta_{oK}^\star,  \alpha_{oH}^\star$, and $\beta_{oH}^\star$ are not shown.

   \begin{figure}[p]
\centering
\includegraphics[width=0.63\linewidth, height=8.5cm]{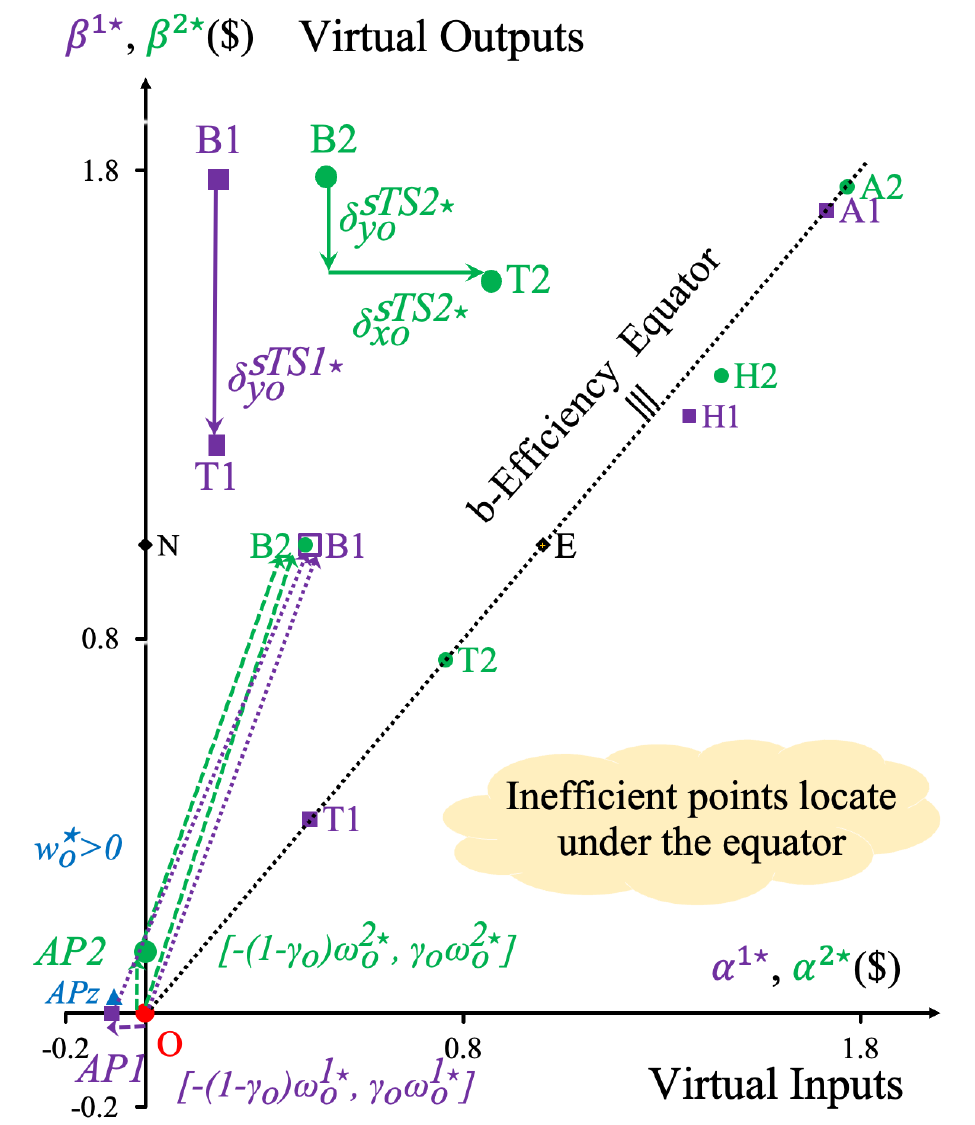} 
\caption{Comparing the solutions of sTS1 and sTS2 models in assessing $DMU_B$.}
\label{fig:5}
\end{figure}

\begin{figure}[p]
\centering
\includegraphics[width=0.66\linewidth, height=8.5cm]{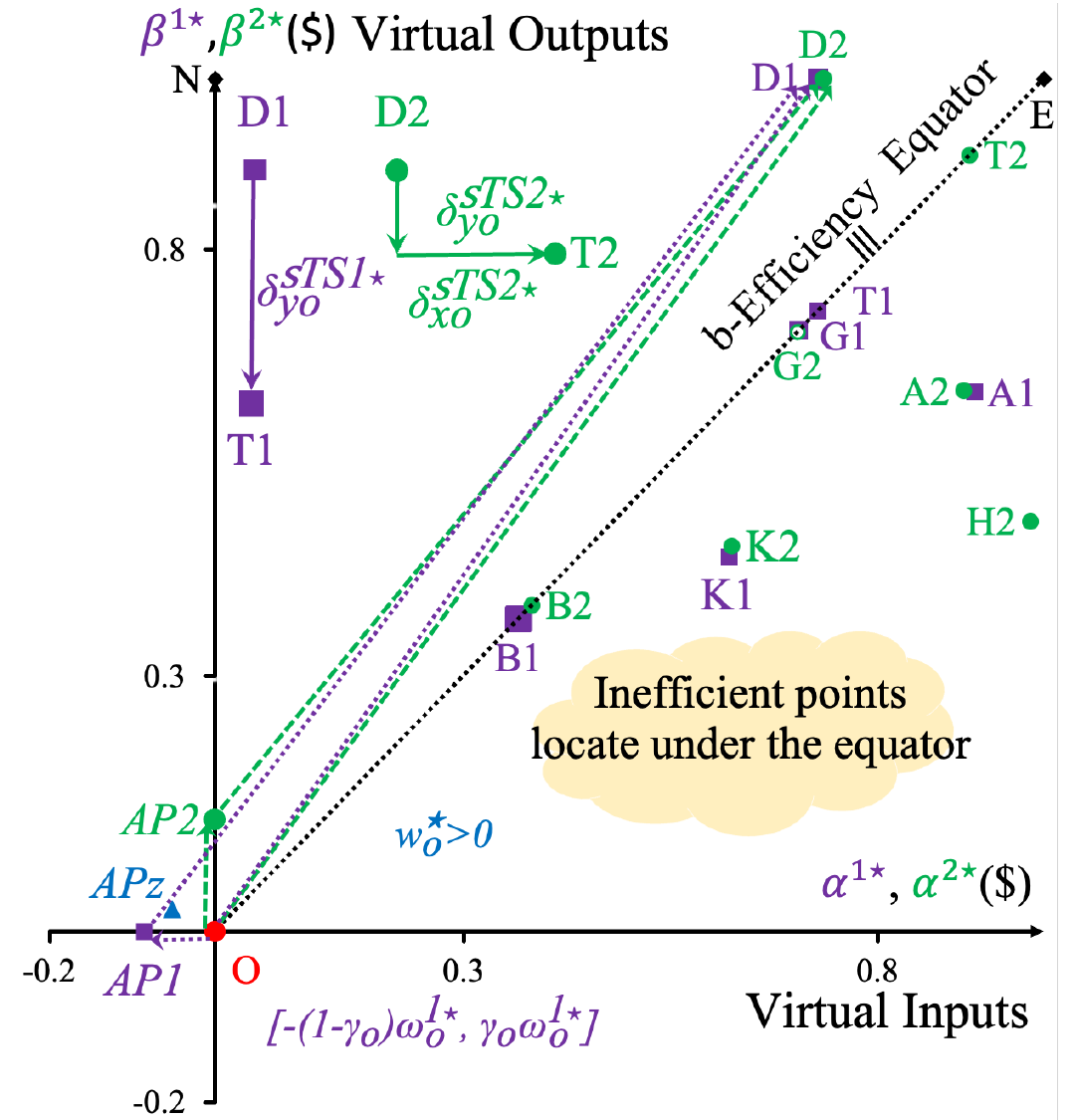}
\caption{Comparing the solutions of sTS1 and sTS2 models in assessing $DMU_D$.}
\label{fig:6}
\end{figure}

\begin{table}[p]
\centering
\setlength{\tabcolsep}{2 pt}
\renewcommand{\arraystretch}{0.91}
\caption{Solutions in Measuring the Super-Efficiencies of $DMU_B$ and $DMU_D$ by Scenario-II.} 
\footnotesize
\begin{tabularx}{\linewidth}{cccrrrrrrrr} \label{table:3} 
 & & $DMU_o$= &$B$ &$B$ &$B$ &$B$ &$D$ &$D$ &$D$ &$D$ \\
\hline 
{Row}&{Sol.}& Unit &sPT-II&sTS1-II&sTSz-II&sTS2-II&sPT-II&sTS1-II&sTSz-II&sTS2-II\\
& & & {Ph.-1} & {Ph.-2}&{Ph.-4}&{Ph.-3}&{Ph.-1}&{Ph.-2}&{Ph.-4}&{Ph.-3}\\
  \hline
R1&$\kappa_o^c $ &  -    & 0.2417 & 0.2417 & 0.3345 & 0.4273 & 1.3411 & 1.3411 & 1.4092 & 1.4774 \\
  &$\tau_o $ &   \$    & 0.5   & 0.5   & 0.4823 & 0.4591 & 0.7709 & 0.8037 & 0.7827 & 0.7614 \\
 & $\Delta_o,\delta_o$& \$ & 0.5855 & 0.5855 & 0.5964 & 0.5979 & 0.2611 & 0.2722 & 0.2688 & 0.2652 \\
 \hline
  R2 &$v_{1o} $  & \$/ton & 0     & 0     & 0     & 0     & 0.3889 & 0.4230 & 0.4119 & 0.4007 \\
   &$v_{2o}$ &   \$/hr    & 0.0143 & 0.0172 & 0.0166 & 0.0158 & 0     & 0     & 0     & 0 \\
   &$u_{1o}$ &    \$/ $m^3$    & 0.0009 & 0.0009 & 0.0009 & 0.0008 & 0.0003 & 0.0003 & 0.0003 & 0.0003 \\
  &$u_{2o}$ &   \$/\%    & 0.0056 & 0.0056 & 0.0054 & 0.0052 & 0.0024 & 0.0020 & 0.0020 & 0.0019 \\
&$w_o$    &    \$   & 0.0000 & 0.3538 & 0.3413 & 0.3249 & 0     & 0.0566 & 0.0551 & 0.0536 \\
          \hline
   R3&$q_{1o}$ &   -    & 0     & 0     & 0     & 0     & 0     & 0.0000 & 0.1156 & 0.2313 \\
   &$q_{2o}$  &   -    & 0     & 0     & 0.3840 & 0.7681 & 0     & 0     & 0     & 0 \\
   &$p_{1o} $   &   -    & 0.4344 & 0.4344 & 0.2172 & 0.0000 & 0.3387 & 0.3387 & 0.2278 & 0.1170 \\
   &$p_{2o} $  &   -    & 0.7366 & 0.7366 & 0.6355 & 0.5343 & 0     & 0     & 0     & 0 \\
&$\pi_{oA}$  &   -    & 0.2417 & 0.2417 & 0.3345 & 0.4273 & 0     & 0     & 0     & 0 \\
&$\pi_{oB}$ &   -    & 0     & 0     & 0     & 0     & 0.6424 & 0.6424 & 0.5211 & 0.3997 \\
&$\pi_{oG}$   &    -   & 0     & 0     & 0     & 0     & 0.6986 & 0.6986 & 0.8881 & 1.0776 \\
\hline
R4 &$\widehat{x}_{1o}$  &   ton    & 1     & 1     & 1     & 51.273 & 1.9   & 1.9   & 2.119 & 2.339 \\
 &$\widehat{x}_{2o}$    &   hr    & 29    & 29    & 40.137 & 567   & 281   & 281   & 281   & 281 \\
& $\widehat{y}_{1o}$   &   $m^3$    & 320.7 & 320.7 & 443.8 & 41.4  & 1617.6 & 1617.6 & 1888.8 & 2159.9 \\
 & $\widehat{y}_{2o}$ &    \%   & 23.442 & 23.442 & 32.444 & 0.5897 & 97    & 97    & 97    & 97 \\
   & $\gamma_o$ &    -   & 0     & 0     & 0.3105 & 0.1388 & 0     & 0.0000 & 0.3367 & 0.6642 \\
  & $\omega_o^c$ &   \$ & 0     & 0.0855 & 0.1142 & 0     & 0     & 0.0759 & 0.0776 & 0.0792 \\
  \hline
R5 &$x_{1o}v_{1o}$&  \$     & 0     & 0     & 0     & 0.4591 & 0.7389 & 0.8037 & 0.7827 & 0.7614 \\
&$x_{2o}v_{2o}$ & \$   &  0.4145 & 0.5   & 0.4823 & 0.4591 & 0     & 0     & 0     & 0 \\
   &$y_{1o}u_{1o}$  &  \$     & 0.5   & 0.5   & 0.4823 & 0.4591 & 0.7709 & 0.8037 & 0.7827 & 0.7614 \\
 &$y_{2o}u_{2o}$ &   \$    & 0.5   & 0.5   & 0.4823 & 0.4021 & 0.2291 & 0.1963 & 0.1912 & 0.1860 \\
\hline
 R6&$\widehat\alpha_o,\widehat\alpha_{o}^c$  &    \$   & 0.4145 & 0.4145 & 0.5888 & 0.7547 & 0.7389 & 0.7278 & 0.8217 & 0.9109 \\
 &$\widehat\beta_o,\widehat\beta_{o}^c$  &    \$   & 0.4145 & 0.4145 & 0.5888 & 0.7547 & 0.7389 & 0.7278 & 0.8217 & 0.9109 \\
 \hline
 R7   & $\Xi_o$  &   -    & 0.4145 & 0.4145 & 0.4036 & 2.4868 & 0.7389 & 0.7278 & 0.7312 & 0.7348 \\
 & $E_o$ &   -    & 2.4126 & 2.4126 & 2.4779 & 2.1621 & 1.3533 & 1.3740 & 1.3676 & 1.3609 \\
 \hline
R8&$\alpha_{oA},\alpha_{oA}^c$ &   \$    & 1.7151 & 1.7151 & 1.7603 & 1.7663 & 0.8945 & 0.9163 & 0.9109 & 0.9037 \\
&$\alpha_{oB},\alpha_{oB}^c$ &   \$    & 0.4145 & 0.4145 & 0.4036 & 0.4021 & 0.3889 & 0.3664 & 0.3754 & 0.3827 \\
&$\alpha_{oD},\alpha_{oD}^c$&    \$   & 4.0163 & 4.4910 & 4.4378 & 3.8242 & 0.7389 & 0.7278 & 0.7312 & 0.7348 \\
&$\alpha_{oG},\alpha_{oG}^c$&     \$  & 3.5732 & 3.9565 & 3.9222 & 1.4497 & 0.7000 & 0.7048 & 0.7049 & 0.7033 \\
  \hline
 R9& $\beta_{oA},\beta_{oA}^c$ &    \$   & 1.7151 & 1.7151 & 1.7603 & 1.7663  & 0.6473 & 0.6323 & 0.6343 & 0.6347 \\
 & $\beta_{oB},\beta_{oB}^c$  &    \$   & 1     & 1     & 1     & 1 & 0.3889 & 0.3664 & 0.3754 & 0.3827 \\
& $\beta_{oD},\beta_{oD}^c$  &   \$    & 2.702 & 2.702 & 2.712 & 1.938 & 1     & 1     & 1     & 1 \\
& $\beta_{oG},\beta_{oG}^c$   &   \$    & 1.902 & 1.902 & 1.941 & 1.362 & 0.7000 & 0.7048 & 0.7049 & 0.7033 \\
          \hline
R10&AP(x)&    \$   & 0     & -0.0855 & -0.0787 & 0.0000 & 0     & -0.076 & -0.051 & -0.027 \\
&AP(y)    &  \$     & 0     & 0     & 0.0354 & 0.1313 & 0     & 0     & 0.026 & 0.053 \\ \hline
\end{tabularx}
\end{table}
\clearpage
 \indent To illustrate the fourth phase, we assume the final SIC scalars are the mid-point of the interval of $\kappa_o^1$ and $\kappa_o^2$. The anchor point of the sTSz model that contains the SIC scalar  $\kappa_o^z$, APz is located in the region of points O, AP1, and AP2 because $DMU_o$ is excluded from Eq. (\ref{Eq:59}). In Figure \ref{fig:5}, AP1 and AP2 are in the second  quadrant because $\omega_o^c>\$0$. When $\omega_o^c>\$0$, the anchor point will be located in the fourth quadrant. R6 in Table \ref{table:3} listed the locations of $DMU_o$ and its projection points on the b-b-Efficiency Equator. $DMU_o$ is located in the first quadrant and above the b-Efficiency Equator. $DMU_o$ projects on the b-Efficiency Equator, the T1 and T2 in Figures \ref{fig:5} and \ref{fig:6}. The reference peers of $DMU_B$, A1, and A2 are on the b-Efficiency Equator, while the remaining inefficient DMUs lie under it. In evaluating $DMU_D$, $DMU_B$, and $DMU_G$, they are located on the b-Efficiency Equator.
 
  \indent In Table \ref{table:3}, $DMU_B$ has higher super-efficiencies (see R7), which reduced the outputs will have an efficiency score equal to 1. Regarding EA, $DMU_B$ has strong outputs and can be selected as the best alternative to the MCDM problem.
    \section {Discussions} \label{sec:6}
 Our EA method effectively estimates the necessary and achievable adjustments for inputs and outputs. Since the b-VGA-EA method is iterative, it can be refined over time as new data is incorporated or the operational context of the DMUs evolves. Our MCDM method assists decision-makers in efficiently selecting the best DMU, contributing significantly to the decision-making community.
 
  \indent The VGA models systematically measured the unified goal price, $\tau_o^\star$ of inputs and outputs  $DMU_o$ in two steps to ensure the objective virtual gap within the range of $(\$0, \$1)$ so that the dimensionless relative efficiencies and relative super-efficiency could be obtained precisely.  

  \indent The current b-VGA-EA method could be modified to the w-VGA-EA method to evaluate each DMU in the light of worst practice. In this context, Scenario I and Scenario II are redefined as Scenario III and Scenario IV, respectively. 
 
\indent  The PT, bTSc, sPT, and sTSc models are modified to wPT, wTSc, hPT, and hTSc models. The wPT and wTSc models estimate the worst-relative efficiency of $DMU_o$, where $E_o^\star>1$. By adjusting the input and output ratios, the worst-relative efficiency will deteriorate from $E_o^\star>1$ to $E_o^\star=1$. 

 \indent The hPT and hTSc models estimate the hypo-relative efficiency of $DMU_o$, where $E_o^\star<1$. Here, adjusting the input and output ratios improves the hypo-relative efficiency from $E_o^\star<1$ to $E_o^\star=1$. In a similar fashion to the b-VGA-MCDM method, the w-VGA-MCDM method can be used to select the worst DMU through a process involving three analogous stages.

\indent The pair of b-VGA-EA and w-VGA-EA methods are both necessary and essential for assessing each DMU. We plan to publish the w-VGA-EA method alongside the w-VGA-MCDM method to provide comprehensive evaluations.

\indent In practical decision-making, it may be necessary to restrict the adjustment ratios. The b-VGA-EA and w-VGA-EA methods can be extended to create several r-VGA-EA methods. For example, by omitting $ \sum_{\forall i\in I} q_{io} \tau_o$ and $ \sum_{\forall r\in R}p_{ro} \tau_o$ in the objective functions of the eight b-VGA and w-VGA models, one can obtain output-oriented and input-oriented measurements. In non-stationary scenarios, where 
$X_a$ and $Y_b$ are restricted from being adjusted, the objective functions of the eight VGA models can be modified to:
 $ \delta_o^{TAP\star} = min \sum_{\forall i\in I\neq a} q_{io} \tau_o + \sum_{\forall r\in R\neq b}p_{ro}\tau_o$.
\indent Typically, the adjustment ratios of inputs ($q_o$) range between 0 and 1, while outputs ($p_o$) may exceed 1. The specific requirement is possible to be met by imposing constraints such as 
($p_o \leq 1$) or define upper and lower bounds:
$\underline{q}_o\le q_o \le \bar{q}_o$ and $\underline{p}_o\le p_o \le \bar{p}_o$. The r-VGA-MCDM methods could be generated to meet the practical conditions.

  \indent The practical issues addressed in the DEA literature\cite{10}, such as the Malmquist Index, network DEA, free-disposal hull, input/output sharing, and input/output selection, can be considered solvable using VGA models. These models can also be adapted to handle datasets containing zeros and negatives\cite{11}. Additionally, VGA models are helpful for datasets with bounded data, ordinal data, and ratio-bounded data\cite{32}.

\indent Further research could explore the development of integrated models that combine the strengths of both MCDM and EA, tailored to specific sectors such as healthcare, energy, and public services. Methodological advancements in addressing uncertainty, adapting to dynamic conditions, and incorporating sustainability dimensions are promising avenues to enhance the applicability and impact of these decision-support tools. In conclusion, the synergy between MCDM methods and EA offers significant potential for advancing sophisticated and practical decision-making frameworks, ultimately improving system performance and supporting the achievement of strategic goals across various contexts.

\indent The TAP program of a VGA model uses a standard form of linear programming with \( (n + m + s) \) columns and \( (m + s + 1) \) rows. This problem can be solved efficiently with free linear programming software on a notebook computer. A software procedure can be created to execute the EA and MCDM methods, ensuring quick execution, including the preparation of EA and MCDM reports and 2D graphical visualizations.
\section* {Acknowledgment} 
We thank the referees for their valuable comments. This research received no specific grant from public, commercial, or not-for-profit funding agencies. We have no competing interests to declare. We confirm that the data supporting the findings of this study are available within the article.

\end{document}